\documentstyle[12pt]{article}
\newcommand{\rbox}[1]{{\rm{#1}}}        
\newcommand{\bbox}[1]{{\bf{#1}}}
\newcommand{\bsym}[1]{{\bf{#1}}}
\newcommand{\rmt}{\rm}
\newcommand{\itt}{\it}
\newcommand{\bft}{\bf}
\def\restriction{\mathbin{\hspace{0.1ex}|\hspace{0.1ex}}}
\def\subsetneqq{\mathbin{\hspace{-0.8mm}
\raisebox{-0.8ex}{
$\msur
\stackrel
{\protect\textstyle\subset}
{\scriptstyle\not=}
\msur$
}\hspace{-0.8mm} }}
\newcommand{\frak}[1]{{\sf{#1}}}
\newcommand{\pri}{{{\bbb\frak{pr}\bbb}_1\hspace{1pt}}}
\newcommand{\prii}{{{\bbb\frak{pr}\bbb}_2\hspace{1pt}}}
\newcommand{\skr}[1]{{\cal #1}}
\newtheorem{theorem}{Theorem}
\newtheorem{assertion}[theorem]{Assertion}
\newtheorem{corollary}[theorem]{Corollary}
\newtheorem{definition}[theorem]{Definition}
\newtheorem{lemma}[theorem]{Lemma}
\newtheorem{proposition}[theorem]{Proposition}
\newtheorem{remark}[theorem]{Remark}
\newcommand{\proof}{\noi{\bft Proof\hspace{2mm} }}
\newcommand{\TF}{\itt}
\newcommand{\bass}{\begin{assertion}\TF\ } 
\newcommand{\eass}{\end{assertion}}
\newcommand{\bcor}{\begin{corollary}\TF\ }
\newcommand{\ecor}{\end{corollary}}
\newcommand{\bdf} {\begin{definition}\rmt\ }
\newcommand{\edf} {\end{definition}} 
\newcommand{\ble} {\begin{lemma}\TF\ }
\newcommand{\ele} {\end{lemma}}
\newcommand{\bpro}{\begin{proposition}\TF\ } 
\newcommand{\epro}{\end{proposition}} 
\newcommand{\brem}{\begin{remark}\rmt\ }
\newcommand{\erem}{\end{remark}} 
\newcommand{\bte} {\begin{theorem}\TF\ }
\newcommand{\ete} {\end{theorem}}
\newcommand{\qed} {\hfill$\msur\Box\msur$} 
\newcommand{\bay}{\begin{array}}
\newcommand{\eay}{\end{array}} 
\newcommand{\bce}{\begin{center}}
\newcommand{\ece}{\end{center}} 
\newcommand{\bde}{\begin{description}}
\newcommand{\ede}{\end{description}}
\newcommand{\ben}{\begin{enumerate}}
\newcommand{\een}{\end{enumerate}}
\newcommand{\bit}{\begin{itemize}}
\newcommand{\eit}{\end{itemize}}

\newcommand{\ZFC} {\bbox{ZFC}}

\newcommand{\cont}[1]{\rbox{Term}_{#1}}
\newcommand{\dom}{\rbox{dom}\,}
\newcommand{\hc}{\rbox{HC}}

\newcommand{\rL}{\rbox{L}}
\newcommand{\ord}{\rbox{Ord}}
\newcommand{\Ord}{\ord}
\newcommand{\od}{\rbox{OD}}
\newcommand{\ran}{\rbox{ran}\,}
\newcommand{\rod}{\hbox{\rmt R\hspace{1pt}-\hspace{1pt}OD}}

\newcommand{\rV}{\rbox{V}}

\newcommand{\sma}[1]{#1\hbox{-}\rbox{SMA}}
\newcommand{\ima}{{\hbox{\hspace{1pt}\rmt ''}}}

\newcommand{\al} {\alpha} 
\newcommand{\ba} {\beta} 
\newcommand{\ga} {\gamma} 

\newcommand{\Da} {\Delta}
\newcommand{\kpa}{\kappa}
\newcommand{\La} {\Lambda} 
\newcommand{\la} {\lambda}

\newcommand{\vpi}{\varphi} 
\newcommand{\om} {\omega} 
\newcommand{\Om} {\Omega} 
\newcommand{\fs}[2]{{\bsym\Sigma}^{#1}_{#2}}
\newcommand{\fp}[2]{{\bsym\Pi}^{#1}_{#2}}
\newcommand{\fd}[2]{{\bsym\Delta}^{#1}_{#2}}
\newcommand{\iSigma}{{\mathchar"7106}}
\newcommand{\is}[2]{\iSigma^{#1}_{#2}}
\newcommand{\iPi}{{\mathchar"7105}}
\newcommand{\ip}[2]{\iPi^{#1}_{#2}}
\newcommand{\iDelta}{{\mathchar"7101}}
\newcommand{\id}[2]{\iDelta^{#1}_{#2}}
\newcommand{\ish}[1]{\is{\rbox{HC}}{#1}}
\newcommand{\iph}[1]{\ip{\rbox{HC}}{#1}}
\newcommand{\idh}[1]{\id{\rbox{HC}}{#1}}
\newcommand{\fdh}[1]{\fd{\rbox{HC}}{#1}}
\newcommand{\bbb}{\hspace{0.5pt}} 
\newcommand{\dF}[1]{{\mathord{{\rm I}\hspace{-2.5pt}{\rm F}}}_{#1}}
\newcommand{\dP}{\mathord{{\rm I}\hspace{-2.5pt}{\rm P}}}
\newcommand{\dX}{\mathord{{\rm X}\hspace{-7pt}{\rm X}}}
\newcommand{\skrsp}{\hspace{0pt}}
\newcommand{\skri}[1]{{\skrsp \skr{#1}\skrsp}}
\newcommand{\cD}{{\skri D}}   
\newcommand{\cH}{H^\ast}      
\newcommand{\cC}[1]{{\cal P}_{#1}}
\newcommand{\cN}{{\skri N}}
\newcommand{\cP}{{\skri P}}
\newcommand{\oP}{{\skri P}^{\hbox{\tiny\rm OD}}}
\newcommand{\cX}{{\skri X}}
\newcommand{\emps}{\emptyset}
\newcommand{\sq}  {\subseteq}

\newcommand{\sneq}{\subsetneqq}
\newcommand{\cj}  {\;\,\&\;\,}

\newcommand{\lra} {\longrightarrow} 
\newcommand{\llra}{\longleftrightarrow} 
\newcommand{\map}{\,\longmapsto\,} 
\newcommand{\res} {{\mathbin{\hspace{0.15ex}\restriction\hspace{0.15ex}}}}
\newcommand{\we}  {{\mathbin{\hspace{0pt}^\wedge}}}
\newcommand{\<}{\leq}
\def\>{\geq}
\newcommand{\ti}{\times}
\newcommand{\dm}{$$}
\newcommand{\ang} [1]{\langle #1\rangle}
\newcommand{\ans} [1]{\{\hspace{0.2mm}#1\hspace{0.2mm}\}}
\newcommand{\mins}{\hspace{-1pt}-\hspace{-1pt}}
\newcommand{\dd}[1]{$\kern-0.7mm{#1}\kern-1mm$-}
\newcommand{\noi} {\noindent}
\newcommand{\vom} {\vspace{1mm}}
\newcommand{\its} {\vspace{-1mm}}
\newcommand{\itla}[1]{\item\label{#1}}
\def\top{{\cal T}}
\def\mattype{\mathbin}
\newcommand{\relfont}[1]{{\sf #1}}
\def\E  {\mattype{\relfont{E}}}

\def\Eo {\mattype{\relfont{E}_0}}
\def\nE {\mattype{\not{\hspace{-2pt}\relfont{E}}}}
\def\nEo{\mattype{\not{\hspace{-2pt}\relfont{E}_0}}}
\def\oE {\mattype{\overline{\relfont{E}}}}
\def\noE{\mattype{\not{\hspace{-2pt}\overline{\relfont{E}}}}}
\def\I{}
\def\J{}
\def\R{}
\def\C{}
\renewcommand{\C}[1]{{\relfont{C}}_{#1}}
\renewcommand{\I}[1]{\mattype{\relfont{Q}_{#1} }}
\newcommand{\Ip}[1] {\mattype{\relfont{Q}'_{#1}}}
\renewcommand{\J}[1]{\mattype{\relfont{R}_{#1} }}
\renewcommand{\R}[1]{\mattype{\relfont{R}_{#1} }}
\def\X{}
\renewcommand{\X}[1]{{\rbox{\sf X}_{#1}}}

\def\t{}
\renewcommand{\t}{{t}}
\newcommand{\tk}[1]{{\t({#1})}}
\newcommand{\ovu}{{\hat u}}
\newcommand{\ovv}{{\hat v}}

\newcommand{\uh}{{\underline{h}}}
\newcommand{\unf}{{\hat{f}}}
\newcommand{\ung}{{\hat{g}}}
\newcommand{\col}[1]{{#1}^{<\om}}
\newcommand{\msur}{\hspace{-1\mathsurround}}

\begin{document}

\normalsize

\title{On a Glimm -- Effros dichotomy and an Ulm--type 
classification in Solovay model}

\author{Vladimir Kanovei
\thanks{Moscow Transport Engineering Institute}
\thanks{{\tt kanovei@math.uni-wuppertal.de} \ and \ 
{\tt kanovei@nw.math.msu.su}
}
\thanks{Partially supported by AMS grant} 
}
\date{01 August 1995} 
\maketitle
\normalsize

\vfill

\begin{abstract}\vspace{2mm} 
\noi
We prove that in Solovay model every $\od$ equivalence $\E$ on 
reals either admits an $\od$ reduction to the equality on the set 
of all countable (of length $<\om_1$) binary sequences, or 
continuously embeds $\Eo,$ the Vitali equivalence.\vspace{1mm}

If $\E$ is a $\is11$ (resp. $\is12$) relation then the reduction 
in the ``either'' part can be chosen in the class of all 
$\id{}1$ (resp. $\id{}2$) functions.
\vspace{1mm}

The proofs are based on a topology generated by $\od$ sets.
\vspace{10mm}

\bce
{\bft Acknowledgements}
\ece
\noi
The author is in debt to M.~J.~A.~Larijani, the president of IPM 
(Tehran, Iran), for the support. 
The author is pleased to thank S.\ D.\ Friedman, G.\ Hjorth, 
\hbox{A.\ S.\ Kechris,} and A.\ W.\ Miller for useful discussions and 
interesting information on Solovay model and the Glimm -- 
Effros matter.
\end{abstract}

\vfill
\vfill

\newpage

\subsection*{Introduction}


The solution of the continuum problem leaves open a variety of 
related questions. For instance, if one works in descriptive set 
theory one may be interested to know how different uncountable 
cardinals can be presented in the real line. This research line can 
be traced as far in the past as the beginning of the century; indeed 
Lebesgue~\cite{leb} found such a presentation for $\aleph_1,$ the 
least uncountable cardinal.

The construction given in \cite{leb} worth to be briefly reviewed. 
One can associate, in an effective way, a set of rationals $Q_x$ 
with each real $x$ so that every set $q$ of rationals has the form 
$Q_x$ for some (perhaps, not unique) $x.$ Let, for a countable 
ordinal $\al$,
\dm
X_\al=\ans{x:Q_x\,\hbox{ is welordered as a set of rationals 
and has the order type }\,\al}\,.
\dm
Then the sets $X_\al,\;\al<\om_1,$ are nonempty and pairwise 
disjoint; therefore we present $\al_1$ in the reals as the sequence 
of the sets $X_\al$. 

This reasoning is a particular case of a much more general 
construction. 

Let $\E$ be an equivalence relation on the reals. Let $\kpa$ be 
the cardinal of the set of all \dd\E equivalence classes; 
then $\kpa\<2^{\aleph_0}.$ One may think  
that the partition of the real line on the \dd\E equivalence 
classes presents the cardinal $\kpa$ in the reals. 

For instance, in the Lebesgue's example, the equivalence can 
be defined as follows: ${x\hspace{2pt}\rL\hspace{2pt} y}$ iff 
either $(1)$ both $Q_x$ and 
$Q_y$ are wellordered and have the same order type, or $(2)$ 
both $Q_x$ and $Q_y$ are {\it not\/} wellordered. The 
\dd\rL equivalence 
classes are the sets $X_\al,\linebreak[3]\;\al<\om_1,$ plus one 
more ``default'' class of all reals $x$ such that $Q_x$ is not 
wellordered.

Of course, one can present {\it every\/} cardinal 
$\kpa\<2^{\aleph_0}$ this way by a suitable equivalence. But the 
problem becomes much more difficult when one works in descriptive 
set theory and looks for an equivalence of a certain ``effective'' 
type. (Notice that the Lebesgue equivalence $\rL$ is a $\is11$ 
relation.) 

This leads us to the following question: given an equivalence 
relation $\E$ on reals, how many equivalence classes it has ? 

The relevant question is then how to ``count'' the classes. 
Generally speaking, counting is a numbering of the given set of 
mathematical objects by mathematical objects of another type, 
usually more primitive in some sense. In particular, the obvious 
idea is to use ordinals (for instance natural numbers) to count 
the equivalence classes. This works well as long as one is not 
interested in the ``effectivity'' of the counting. Otherwise we 
face problems even with very simple relations. (Consider the 
equality as an equivalence relation. Then one cannot define in 
$\ZFC$ an ``effective'' in any reasonable sense counting of the 
equivalense classes, alias reals, by ordinals.) 

The other natural possibility is to use {\it sets of ordinals\/} 
(for instance reals) to count the equivalence classes.~\footnote
{\rmt\ There are known many mathematical examples, in probability 
and the measure theory, based on this type of enumeration of 
the equivalence classes, see Harrington, Kechris, and Louveau 
\cite{hkl}.} Note that the next step, that is, counting by 
{\it sets of sets\/} of ordinals, would be silly because the 
classes themselves are of this type.

\bdf
[\hspace{1pt}Informal\hspace{1pt}]\\[1pt] 
An equivalence relation is 
{\it discrete\/} iff it admits an ``effective'' enumeration of 
the equivalence classes by ordinals. An equivalence relation is 
{\it smooth\/} iff it admits an ``effective'' enumeration of the 
equivalence classes by {\it sets of\/} ordinals.\qed
\edf
Of course the definition has a definite meaning only provided 
one makes clear the meaning of the ``effectivity''. However in 
any reasonable case one has the following two counterexamples:
\vspace{3mm}

\noi
{\sl Example 1\/}. The equality on a perfect set of reals is 
{\it not\/} discrete.\vspace{2mm}

\noi
{\sl Example 2\/}. The Vitali equivalence relation is 
{\it not\/} smooth.\vspace{3mm}

\noi 
({\it Not\/} here means that one cannot prove in $\ZFC$ the 
existence of the required enumerations among the 
real--ordinal definable functions. However different additional 
axioms, e. g. the axiom of constructibility, make each 
equivalence discrete in certain sense.) 

At the first look, there should be plenty of other counterexamples. 
However, in certain particular but quite representative cases one 
can prove a {\it dichotomy theorem\/} which says that an 
equivalence relation is not discrete (resp. smooth) iff it contains 
Example 1 (resp. Example 2). This is also the topic of this article, 
but to proceed with the reasoning we need to be more exact. 

\subsubsection*{Notation}

Let us review the basic notation of Harrington, Kechris, and 
Louveau \cite{hkl}. See \cite{hkl} or \cite{hk} for a more 
substantial review with details and explanations.

Let $\E$ and $\E'$ be equivalence relations on resp. sets $X,\,X'$. 

A function $U:X\;\longmapsto\;X'$ is a {\it reduction\/} of $\E$ to 
$\E'$ iff $x\E y\;\llra\;U(x)\E' U(y)$ holds for all $x,\,y\in X.$ 
An {\it enumeration of the\/ \dd\E equivalence classes\/} 
(by elements of $X'$) is a reduction of $\E$ to the equality on $X'.$ 
(Here, it is not assumed that all of elements of $X'$ are involved.) 
 
A $1-1$ reduction is called an {\it embedding\/}. $\E'$ 
{\it continuously embeds\/} $\E$ iff there exists a continuous 
embedding $\E$ to $\E'.$ In the case when $X$ is the {\it Cantor 
set\/} $\cD=2^\om$ (with the product topology), $\E'$ 
continuously embeds $\E$ if and only if there exists a perfect set 
$P\sq X$ such that $\ang{P\hspace{1pt};\hspace{1pt}\E\res P}$ is 
homeomorphic to $\ang{X'\hspace{1pt};\hspace{1pt}\E'}.$ In other 
words, embedding $\E$ continuously means in this case that $\E'$ 
contains a homeomorphic copy of $\E$. 

In particular $\E'$ continuously embeds the equality on $\cD$ 
iff there exists a perfect set of \dd{\E'}inequivalent points.

Finally, let $\Eo$ denote the {\it Vitali equivalence\/} on 
$\cD=2^\om,$ defined as follows: ${x\Eo y}$ iff $x(n)=y(n)$ for 
almost all (i.\ e. all but finite) $n\in\om.$

\subsubsection*{The main theorem}

This paper intends to complete the diagram of the following 
three classical theorems on equivalence relations. 
\bde
\item [$\hbox{\sl Borel -- 1\/}.$] Each Borel equivalence on 
reals, either has countably many equivalence classes or admits 
a perfect set of pairwise inequivalent points. (Silver~\cite{si}, 
in fact for \dd{\fp11}relations.)\its

\item [$\hbox{\sl Borel -- 2\/}.$] Each Borel equivalence relation 
on reals, either admits a Borel enumeration of the equivalence 
classes by reals~\footnote
{\rmt\ That is, admits a Borel reduction to the equality on reals. 
Such an equivalence is called {\em smooth\/} in the notation of 
\cite{hkl}.}   
, or continuously embeds the Vitali equivalence $\Eo.$ 
(The Glimm -- Effros dichotomy 
theorem of Harrington, Kechris, and Louveau~\cite{hkl}.)\its
%
\item [$\hbox{\sl Solovay model -- 1}.$] In Solovay 
model~\footnote
{\rmt\ By Solovay model we mean a generic extension $\rL[G]$ of 
$\rL,$ 
the class of all constructible sets, by a generic over $\rL$ 
subset of a certain notion of forcing $\cP^\Om\in\rL$ which 
provides a collapse of all cardinals in $\rL,$ smaller than a 
fixed inaccessible cardinal $\Om,$ to $\om,$ see 
Solovay~\cite{sol}. In this model, all projective sets are 
Lebesgue measurable.}
, each 
$\rod$ (real--ordinal definable) equivalence on $\cN$ either 
has $\<\aleph_1$ equivalence classes and admits a $\rod$ 
enumeration of them, or admits a perfect set of 
pairwise inequivalent points. 
(Stern~\cite{ste}.)
\ede
Thus the results $\hbox{\sl Borel -- 1\/}$ and 
$\hbox{\sl Solovay model -- 1}$ say (informally) that an 
equivalence relation either is discrete or contains a 
continuous copy of Example~1 above. Similarly 
$\hbox{\sl Borel -- 2\/}$ says that an 
equivalence relation either is smooth or contains a 
continuous copy of Example~2 above.
\bte
\label{main}
{\rmt [{\sl\hspace{1pt}Solovay model -- 2\/\hspace{1pt}}]}\\[1pt]
The following is true in Solovay model. Assume that\/ $\E$ is an\/ 
$\rod$ equivalence on\/ $\cN.$ Then one and 
only one of the following two statements holds$:$\its
\ben
\def\theenumi{\hskip2pt{\rmt(\Roman{enumi})}\hskip2pt}
\def\labelenumi{\theenumi}
\itla{1}\msur 
$\E$ admits a\/ $\rod$ enumeration of the equivalence classes by 
elements of\/ $2^{<\om_1}.$~\footnote
{\rmt\ $2^{<\om_1}=\bigcup_{\al<\om_1}2^\al$ denotes the set of all 
countable (of any length $<\om_1$) binary sequences.} 

If moreover $\E$ is a\/ $\fs11$ {\rmt({\it resp. $\fs12$})} 
equivalence then the enumeration exists in the 
class\/ $\fdh1$ {\rmt({\it resp. $\fdh2$})}~\footnote
{\rmt\ Here by $\fdh n$ we denote the class of all subsets of 
$\hc$ (the family of all hereditarily countable sets) which are 
$\id{}n$ in $\hc$ by formulas which may contain arbitrary 
{\it reals\/} as parameters.}$;$\its

\itla{2}\msur 
$\E$ continuously embeds $\Eo$.\its
\een
\ete
\nopagebreak
This is the main result of this paper.\vom\pagebreak[3]

{\it Remark 1\/}. Hjorth~\cite{h-det} obtained a similar theorem 
in a strong determinacy hypothesis (${\bf AD}$ holds in 
$\rL[\hbox{reals}]$), yet with a weaker part \ref{1}: an $\od$ 
reduction to the equality\pagebreak[3] 
on a set $2^\kappa,$ $\kappa\in\Ord$.
\vom

{\it Remark 2\/}. The statements \ref{1} and \ref{2} are 
incompatible. Indeed otherwise there would exist an $\rod$ 
function $U:\cD\,\lra\,2^{<\om}$ which reduces $\Eo$ to the 
equality on $2^{<\om}.$ Let $U$ be $\od[z],$ $z\in\cD.$ Then for 
each $p\in\ran U\;\,(\sq 2^{<\om}),$ $F(p)=U^{-1}(p)$ is a 
\dd{\Eo}equivalence class, a countable $\od[p,z]$ subset of 
$\cD.$ In Solovay model, this implies $F(p)\sq \rL[z,p]$ for 
all $p.$ We obtain an $\od[z]$ choice function 
${g:\ran U\;\lra\;\cD}$ such that $g(p)\in F(p)$ for all $p.$ 
Then $\ran p$ is an $\rod$ selector for $\Eo,$ hence a 
non\-measurable $\rod$ set, contradiction with the known 
properties of the Solovay model.\vom

{\it Remark 3\/}. $2^{<\om_1}$ cannot be replaced in 
Theorem~\ref{main} by an essentially smaller set. To see this  
consider the equivalence $\R{}$ on $\cN^2$ defined as follows: 
$\ang{z,x}\R{}\ang{z',x'}$ iff\its
\bit
\item[--] {\it either\/} $z$ and $z'$ code the same 
countable ordinal and $x$ and $x'$ code, in the sense of 
$z$ and $z'$ respectively, the same subset of the ordinal, \its

\item[--] {\it or\/} both $z$ and $z'$ do not code an ordinal.\its
\eit
The relation $\R{}$ admits an $\od$ reduction \underline{on}to 
$\Da(2^{<\om_1}),$ the equality on $2^{<\om_1},$ therefore does 
not embed $\Eo$ continuously in Solovay model (see Remark 1). It 
follows that any set $W$ such that $\R{}$ 
admits a $\rod$ reduction to $\Da(W)$ has a subset $W'\sq W$ 
which is in $1-1$ $\rod$ correspondence with $2^{<\om_1}.$ In 
particular, the continuum $\cN$ does not satisfy this condition 
in Solovay model. (Indeed $2^{<\om_1}$ has $\rod$ subsets of 
cardinality exactly $\aleph_1$ while $\cN$ does not have those in 
Solovay model.)\vom 

{\it Remark 4\/}. Even in the case of $\is11$ equivalence 
relations $2^{<\om_1}$ cannot be replaced by $\cN$ in \ref{1}. 
Indeed the $\is11$ equivalence $x\E y$ iff 
either $x,\,y\in\cN$ code the same (countable) ordinal 
or both $x$ and $y$ do not code an ordinal 
(Example 6.1 in Hjorth and Kechris~\cite{hk}) neither admits a 
$\fd12$ reduction to $\Da(\cN)$ nor embeds $\Eo$ via a $\fd12$ 
function in $\ZFC$ plus 
$\forall\,x\in\cN\,(\om_1^{\rL[x]}<\om_1).$ (In Solovay model, 
$\fd12$ can be strengthened to $\rod$.) This shows that the 
Glimm -- Effros theorem of Harrington, Kechris, and Louveau 
\cite{hkl} (theorem $\hbox{\sl Borel -- 2\/}$ above) cannot be 
expanded from Borel to $\fs11$ relations.\vom

{\it Remark 5\/}. On the other hand, $\fs11$ equivalence relations 
tend to satisfy a looser {\em Ulm--type\/} dichotomy.~\footnote
{\ The notion introduced in~\cite{hk}. Hjorth and Kechris 
refer to certain classification results in algebra, i.e. the Ulm 
classification of countable abelian \dd pgroups.} 
In particular, Hjorth and Kechris~\cite{hk} proved that every 
$\fs11$ equivalence {\it with Borel classes\/} either admits a 
$\fd{}1$ reduction to $\Da(2^{<\om_1}),$ the equality on 
$2^{<\om_1},$ or embeds $\Eo$ continuously; furthermore in the 
assumption $\forall\,x\in\cN\,(x^{\#}\hbox{ exists})$ the 
requirement that the \dd\E classes are Borel can be dropped. 

Thus Theorem~\ref{main} proves that the Ulm classification is 
available in the Solovay model. This gives a partial answer to the 
question posed by Hjorth and Kechris in~\cite{hk}.~\footnote
{\rmt\ ``Is $\forall\,x\in\cN\;(x^{\#}\,\hbox{ exists})$ needed 
to prove that'' \ref{1} (with a $\fdh1$ reduction) or \ref{2} 
hold for $\fs11$ relations; item 3) in Section 7: 
{\sl Open problems\/} 
in a preprint version of~\cite{hk}. Since the sharps hypothesis 
fails in Solovay model, we observe that the answer is: {\em not}.}

It would be interesting to get the Ulm classification 
for $\fs11$ relations in $\ZFC$.~\footnote
{\rmt\ The author~\cite{k-s11} proved the result assuming that 
each real belongs to a set--generic extension of $\rL$.}



\subsubsection*{A brief description of the exposition}

\bde


\item[{\itt Section \ref{appr}$:$}]
We outline how the proof of Theorem~\ref{main} will go on. 
A topology $\top$ generated by $\od$ sets in Solovay model (a 
counterpart of the Gandy -- Harrington topology) is introduced. 
Similarly to Harrington, Kechris, and Louveau~\cite{hkl}, we have 
two cases: either the equivalence $\E$ of consideration is closed 
in the topology $\top^2$ or it is not closed. The plan of the 
proof of Theorem~\ref{main} is to demonstrate that the first case 
provides \ref{1} while the second provides \ref{2}. 
We also review some important properties of the Solovay model. 

\item[{\itt Section \ref{clos}$:$}] 
We prove that in the case when $\E=\oE$ the equivalence $\E$ 
satisfies the requirements of Item~\ref{1} of Theorem~\ref{main}. 
The argument for the ``moreover'' part of Item~\ref{1} includes 
the idea of forcing the equivalence of mutually generic reals over 
countable models, due to Hjorth and Kechris~\cite{hk}.

\item[{\itt Section \ref{top}$:$}]
We begin to study the case when the given equivalence is not 
\dd{\top^2}closed in Solovay model. We develop forcing notions 
$\dX$ and $\dP$ associated with $\top$ and $\top^2$ respectively. 
In particular it is demonstrated that the intersection of a 
generic set is nonempty. The set 
$H=\ans{x:[x]_{\E}\sneq [x]_{\oE}},$ nonempty 
as soon as we assume $\E\sneq\oE,$ is considered. 

\item[{\itt Section \ref{or}$:$}]
We accomplish the case when the given relation $\E$ is 
not \dd{\top^2}closed. It is demonstrated that in this case $\E$ 
continuously embeds $\Eo.$ The splitting construction is based on 
the principal idea of Harrington, Kechris, and Louveau \cite{hkl}, 
but the technical realization is quite different since we use 
straightforward forcing arguments rather than Choquet games, which 
makes the construction a little bit more elementary.
\ede

\noi
{\bf Important remark} \\[1mm]
It will be convenient to use the {\em Cantor set} $\cD=2^\om$ 
rather than the Baire space $\cN=\om^\om$ as the principal space 
in this paper. 

\newpage

\subsection{Approach to the main theorem}
\label{somo}
\label{appr}

In this section, we outline the proof of Theorem~\ref{main}, 
the principal theorem of 
the paper. The idea of the proof has a semblance of the proof of 
the ``Borel'' Glimm -- Effros theorem in Harrington, Kechris, and 
Louveau \cite{hkl}; in particular the dichotomy will be determined 
by an answer to the question whether the given relation $\E$ is 
closed in a certain topology on $\cD^2$.

First of all we review the definition and some properties of 
Solovay model.

\subsubsection{The Solovay model}
\label{bnot}

Let $\al$ be an ordinal. Then $\col\al$ is the forcing to collapse 
$\al$ down to $\om.$ 
We let $\cC{<\la}$ be the product of all sets $\col\al,\;\,\al<\la,$ 
with finite support; in other words, $\cC{<\la}$ is the set of all 
functions $p$ defined on finite subsets of $\la$ such that 
$p(\al)\in\col\al$ for each $\al<\la,$ $\al\in\dom p$. 

The forcing notions $\col\al$ and $\cC{<\la}$ are equivalent 
respectively to 
${\cal P}_\al$ and ${\cal P}^\la$ in Solovay~\cite{sol}. 
We set $\cC{\<\la}=\cC{<{\la+1}}$.

(Notice that the definitions of 
$\cC{<\la}$ and $\cC{\<\la}$ are absolute.)

Let $M$ be a transitive model of $\ZFC,$ a set or proper class, 
containing $\Om,$ an inaccessible cardinal in $\rL.$ 
By \dd\Om{\em Solovay extension of $M$} we shall 
understand a generic extension of the form $M[G],$ where 
$G\sq\cC{<\Om}$ is \dd{\cC{<\Om}}generic over $M$. 

\bdf
\label{dfsma}
\dd\Om{\em Solovay model axiom}, $\sma\Om$ in brief, is the 
following hypothesis:

{\it $\Om$ is inaccessible in $\rL$ and the universe $\rV$ is an 
\dd\Om Solovay extension of $\rL$.}\qed
\edf

\subsubsection{The dichotomy}
\label{dichot}

As usual, we shall concentrate on the ``lightface'' case of an 
$\od$ equivalence relation $\E;$ the general case when $\E$ is 
$\od[z]$ for a real $z$ can be carried out similarly. 

Thus let us consider an $\od$ equivalence $\E$ in the assumption 
of $\sma\Om$. 

The relation $\E$ is fixed in the remainder of the proof of 
Theorem~\ref{main}. The hypothesis $\sma\Om$ will be assumed during 
the proof of the theorem, but we 
shall not mind to specify $\sma\Om$ explicitly in all formulations 
of of lemmas etc. 

For any set ${X\sq\cD},$ we put 
$[X]_{\E}=\ans{y:\exists\,x\in X\;(x\E y)},$ the 
{\em\dd\E saturation} of $X$. 

Let $\top$ be the topology generated on a given set $X$ (for 
instance, $X=\cD=2^\om,$ the Cantor set) by all $\od$ 
subsets of $X.$ $\top^2$ is the product of two copies of $\top,$ 
a topology on $\cD^2$.

We define $\oE$ to be the \dd{\top^2}closure of $\E$ in $\cD^2.$ 
Thus ${x\noE y}$ iff there exist $\od$ sets $X$ and $Y$ containing 
resp. $x$ and $y$ and such that ${x'\nE y'}$ for all ${x'\in X,}$ 
${y'\in Y}.$ Obviously $X$ and $Y$ can be chosen as \dd\E 
invariant sets (otherwise take \dd\E saturations of $X$ and $Y$), 
and then $Y$ can be replaced by the complement of $X,$ so that 
\dm
x\oE y\;\;\llra\;\;\forall\,X\;
[\,X\hbox{ is }\od \cj X\hbox{ is \dd\E invariant}\;\,
\lra\;\,(x\in X\;\llra\;y\in X)\,]\,.
\dm
Therefore $\oE$ is an $\od$ equivalence, too. 

We now come to the key point of the dichotomy: either $\E=\oE$ 
or $\E\sneq\oE$. This reduces the ``lightface'' case in 
Theorem~\ref{main} to the following form. 

\bte
\label{mt}
Assume\/ $\sma\Om.$ Let\/ $\E$ be an\/ $\od$ equivalence on\/ 
$\cD.$ Then\its
\ben
\def\theenumi{\hskip2pt{\rmt(\Roman{enumi})}\hskip2pt}
\def\labelenumi{\theenumi}
\itla{1t} If\/ $\E=\oE$ then\/ $\E$ admits an\/ $\od$ enumeration 
of the equivalence classes by elements of \/ $2^{<\om_1}$.

If moreover $\E$ is a\/ $\is11$ {\rmt({\it resp. $\is12$})} 
equivalence then the enumeration exists in the class\/ $\idh1$ 
{\rmt({\it resp. $\idh2$})}$;$\its

\itla{2t} If\/ $\E\sneq\oE$ then\/ $\E$ continuously embeds\/ 
$\Eo$. 
\een
\ete
In the case when the relation $\E$ is $\od[z]$ (resp. 
$\is11[z]\,,\;\,\is12[z]$) for a real $z,$ the $z$ uniformly 
enters the reasoning, not causing any problem. (In particular 
one considers $\top[z],$ the topology generated by $\od[z]$ 
sets, rather than $\top$.) 

We prove part \ref{1t} of the theorem in the next section.  
Part \ref{2t} will be considered in the two following sections. 
The rest of this section presents different properties of the 
Solovay model.

\subsubsection{Weak sets in Solovay model}
\label{weak}

A set $x$ will be called \dd\Om{\em weak over $M$} iff $x$ belongs 
to an \dd{\col\al}generic extension of $M$ for some $\al<\Om$.

\bpro
\label{solM}
Assume\/ $\sma\Om.$ Then\/ $\Om=\om_1.$ Furthermore, suppose 
that\/ $S\sq\Ord$ is\/ \dd\Om weak over\/ $\rL.$ Then\/ \its
\ben
\def\theenumi{{\arabic{enumi}}}
\def\labelenumi{{\rmt\theenumi.}}
\itla{sm1}
$\Om$ is inaccessible in\/ $\rL[S]$ and\/ 
$\rV$ is an\/ \dd\Om Solovay extension of $\rL[S]$.\its

\itla{sm2} 
If\/ $\Phi$ is a sentence containing only sets in\/ $\rL[S]$ as 
parameters then\/ $\La$ {\rm(}the empty function\/{\rm)} decides\/ 
$\Phi$ in the sense of\/ $\cC{<\Om}$ as a forcing notion over 
$\rL[S]$.\its

\itla{sm4} 
If a set\/ $X\sq\rL[S]$ is\/ $\od[S]$ then\/ 
$x\in\rL[S]$.
\een 
\epro
($\od[S]=S$\msur{}--{\it ordinal definable\/}, that is, 
definable by an \dd\in formula containing $S$ and ordinals as 
parameters.)

The proof (a copy of the proof of Theorem 4.1 in Solovay~\cite{sol}) 
is based on several lemmas, including the following crucial lemma:

\ble
\label{44}
{\rmt (Lemma 4.4 in \cite{sol})}\\[1pt]
Let\/ $M$ be a transitive model of\/ $\ZFC,$ $\la\in\Ord\cap M.$ 
Suppose that\/ $M'$ is a\/ \dd{\col\la}generic extension of\/ 
$M$ and\/ $M''$ is a\/ \dd{\col\la}generic extension of\/ $M'.$ 
Let\/ $S\in M',$ $S\sq\Ord.$ Then\/ $M''$ is a\/ 
\dd{\col\la}generic extension of\/ $M[S]$.\qed
\ele
\proof{}of the proposition. 

{\em Item \ref{sm1}\/}. By definition, $S$ belongs to a 
\dd{\col\al}generic extension of $\rL$ for some $\al<\Om.$ Then in 
fact $S\in\rL[x]$ for a real $x.$ It follows (Corollary 
3.4.1 in \cite{sol}) that there exists an ordinal $\la<\Om$ such 
that $S$ belongs to the model $M_\la=\rL[G_{\<\la}]$ for some 
$\la<\Om,$ where $G_{\<\la}=G\cap\cC{\<\la}$. 

Notice that $G_{\<\la}$ is \dd{\cC{\<\la}}generic over $\rL.$ 
Therefore by Lemma 4.3 in Solovay~\cite{sol}, $M'=M_\la$ is a 
\dd{\col\la}generic extension of $\rL$. 

Let us consider the next step $\la+1.$ 
Obviously the model 
$M_{\la+1}=\rL[G_{\<\la+1}]$ is a \dd{\col{(\la+1)}}generic 
extension of $M_\la.$ Since $\col{(\la+1)}$ is order isomorphic to 
the product $\col\la\ti\col{(\la+1)},$~\footnote
{\ Here $\col\la\ti\col{(\la+1)}$ is the set of all pairs 
$\ang{p,q}$ such that $S\in \col\la,$ $q\in\col{(\la+1)},$ and 
$\dom p=\dom q$.}
we conclude that $M_{\la+1}$ is a \dd{\col{(\la+1)}}generic 
extension of a certain \dd{\col\la}generic extension $M''$ of 
$M'=M_\la$. 

Lemma \ref{44} says that $M''$ is a \dd{\col\la}generic extension 
of $\rL[S],$ therefore a \dd{\cC{\<\la}}generic extension of 
$\rL[S]$ as well by Lemma 4.3 in \cite{sol}. 

It follows that $M_{\la+1}$ is a \dd{\cC{\<\la+1}}generic 
extension of $\rL[S].$ 

Finally $M=\rL[G]$ is a \dd{\cC{\>\la+2}}generic extension of 
$M_{\<\la+1}=\rL[G_{\<\la+1}].$ This ends the proof of item 
\ref{sm1} of the proposition.

{\em Items \ref{sm2} and \ref{sm4}\/}. It 
suffices to refer to item \ref{sm1} and apply resp. Lemma 3.5 and 
Corollary 3.5 in \cite{sol} for $\rL[S]$ as the initial model.\qed

\subsubsection{Coding of reals and sets of reals in the model}
\label{coding}


If $G\sq\col\al$ is \dd{\col\al}generic over a transitive model 
$M$ ($M$ is a set or a class) then $f=\bigcup G$ maps $\om$ onto 
$\al,$ so that $\al$ is countable in $M[G]=M[f].$ Functions 
$f:\om\,\lra\,\al$ obtained 
this way will be called \dd{\col\al}{\em generic over\/ $M$.}

We let $\dF\al(M)$ be the set of all \dd{\col\al}generic over $M$ 
functions $f\in \al^\om.$ We put $\dF\al[S]=\dF\al(\rL[S])$ and 
$\dF\al=\dF\al(\rL)=\dF\al[\emps]$.

We recall that $\cD=2^\om$ is the principal descriptive space in 
this research. The following definitions intend to give a useful 
coding system for reals and sets of reals in Solovay model.

Let $\al\in\Ord.$ By $\cont\al$ we denote the set of all 
``terms'' --- indexed 
sets $\t=\ang{\al,\ang{\t_n:n\in\om}}$ such that 
${\t_n\sq\col\al}$ for each $n$. 

We put $\cont{}=\bigcup_{\al<\om_1}\cont\al.$ (Recall that 
$\om_1=\Om$ assuming $\sma\Om$.)  

``Terms'' $\t\in\cont\al$ are used to code functions 
$C:\al^\om\;\lra\;\cD=2^\om;$ namely, for every $f\in\al^\om$ 
we define $x=\C\t(f)\in\cD$ by: $x(n)=1$ iff $f\res m\in \t_n$ 
for some $m$.

Assume that $\t=\ang{\al,\ang{\t_n:n\in\om}}\in\cont\al,$  
$u\in\col\al,$ $M$ arbitrary. We introduce the sets 
$\X{\t u}(M)=\ans{\C\t(f):u\subset f\in\dF\al(M)}$ and 
$\X\t(M)=\X{\t \La}(M)=\C\t\ima\dF\al(M).$ As above, we let 
$\X\t[S]=\X\t(\rL[S])$ and $\X\t=\X\t[\emps]=\X\t(\rL);$ the same 
for $\X{\t u}$.

\bpro
\label{solMb}
Assume\/ $\sma\Om.$ Let $S\sq\Ord$ be\/ \dd\Om weak over\/ 
$\rL.$ Then\its
\ben
\def\theenumi{{\arabic{enumi}}}
\def\labelenumi{{\rmt\theenumi.}}
\itla{sm6} 
If\/ $\al<\Om,$ $F\sq\dF\al[S]$ is\/ $\od[S],$ and\/ $f\in F,$ 
then there exists\/ $m\in\om$ such that each\/ $f'\in\dF\al[S]$ 
satisfying\/ $f'\res m= f\res m$ belongs to\/ $F$.\its

\itla{sm5} 
For each\/ $x\in\cD,$ there exist\/ $\al<\Om=\om_1,$ 
$f\in\dF\al[S],$ and\/ $\t\in\cont\al\cap\rL[S]$ such that\/ 
$x=\C\t(f)$.\its

\itla{xl1} 
Each\/ $\od[S]$ set $X\sq\cD$ is a union of sets of the 
form\/ $\X\t[S],$ where $\t\in\cont\al\cap\rL[S]$ for some\/ 
$\al<\Om=\om_1$.\its

\itla{xl2} 
Suppose that\/ ${\t\in\cont\al\cap\rL[S],\;\;\al<\Om=\om_1,}$ and\/ 
${u\in\col\al}.$ Then every\/ $\od[S]$ set\/ 
$X\sq\linebreak[3]{\X{\t u}[S]}$ is a union of sets of 
the form\/ $\X{\t v}[S],$ where\/ $u\sq v\in\col\al$.
\een
\epro
\proof {\em Item \ref{sm6}\/}. We observe that  
$F=\ans{f'\in\al^\om:\Phi(p,f')}$ for an \dd\in formula $\Phi.$ 
Let $\Psi(p,f')$ denote the formula: ``$\La$ \dd{\cC{<\Om}}forces 
$\Phi(p,f')$ over the universe'', so that 
\dm
F=\ans{f'\in\al^\om:\Psi(p,f')\,\hbox{ is true in }\,\rL[S,f']}.
\dm
by Proposition~\ref{solM} (items \ref{sm1} and \ref{sm2}). 
Therefore, since $f\in F\sq\dF\al[S],$ there exists $m\in\om$ such 
that the restriction $u=f\res m$\hspace{1mm}{} 
(then $u\in\col\al$)\hspace{1mm}{} 
\dd{\col\al}forces $\Psi(p,{\hat f})$ over $\rL[S],$ 
where $\hat f$ is the name of the \dd\al collapsing function. 
The $m$ is as required. 

{\em Item \ref{sm5}\/}. Since the universe is a Solovay extension 
of $\rL[S]$ (Proposition~\ref{solM}), $x$ belongs to an 
\dd{\col\al}generic extension of $\rL[S],$ for some $\al<\Om.$ 
Thus $x\in\rL[S,f]$ where $f\in\dF\al[S].$ Let 
${\hat x}$ be the name of $x.$ We put
$\t_n=\ans{u\in\col\al:u\,\hbox{ forces }\,{\hat x}(n)=1}$.

{\em Item \ref{xl1}\/}. Let $x\in X.$ We use item 2 to get 
$\al<\Om,$ $f\in\dF\al[S],$ and\/ $\t\in\cont\al\cap\rL[S]$ such 
that\/ $x=\C\t(f).$ Then we apply item 1 to the $\od[S]$ set 
\dm
F=\ans{f'\in\dF\al[S]:\C\t(f')\in X}
\dm 
and the given function $f.$ 
This results in a condition $u=f\res m\in\col\al$ ($m\in\om$) 
such that $x\in \X{\t u}[S]\sq X.$ Finally the set 
$\X{\t u}[S]$ is equal to $\X{\t'}[S]$ for some other 
$\t'\in \cont\al\cap\rL[S]$.

{\em Item \ref{xl2}\/}. Similar to the previous item.\qed



\newpage

\subsection{The case of a closed relation}
\label{clos}

In this section, we prove item \ref{1t} of Theorem~\ref{mt}. 
Thus let us suppose $\sma\Om$ and consider an $\od$ equivalence 
relation $\E$ on $\cD$ satisfying $\E=\oE$.

First of all we obtain a characterization for $\oE$.

We recall that $\Om=\om_1$ in the assumption $\sma\Om,$ and 
$\cont{}=\bigcup_{\al<\om_1}\cont\al$. 

Let us fix an $\od$ enumeration 
$\cont{}\cap\rL=\ans{\t(\xi):\xi<\om_1}$ such that each ``term'' 
$\t\in\cont{}\cap\rL$ has uncountably many numbers $\xi,$ and 
$\t(\xi)\in\cont\al$ for some $\al\<\xi$ 
whenever $\xi<\om_1=\Om$. 
%

\ble
\label{har}
Assume\/ $\sma\Om$ and\/ $\E=\oE.$ Let\/ $x,\,y\in\cD.$ Then\/ 
${x\E y}$ is equivalent to each of the following two 
conditions$:$\\[4mm]
$
\bay{rrclll}
(i) & x\in[\X{\t(\xi)}(\rL)]_{\E} & \llra & 
y\in[\X{\t(\xi)}(\rL)]_{\E} & \hbox{ for each} & \xi<\om_1\,;\\[2mm]

(ii) & x\in[\X{\t(\xi)}(\rL_\xi)]_{\E} & \llra & 
y\in[\X{\t(\xi)}(\rL_\xi)]_{\E} & \hbox{ for each} & \xi<\om_1\,.
\eay
$
\ele
\proof $x\E y$ implies both $(i)$ and $(ii)$ because the sets 
$\X{\t(\xi)}(\rL)$ and $\X{\t(\xi)}(\rL_\xi)$ are $\od.$ Let us 
prove the opposite direction. 

Assume that ${x\noE y}.$ There exists an $\od$ set $X$ such that 
$x\in [X]_{\E}$ but $y\not\in [X]_{\E}.$ By 
Proposition~\ref{solMb}, $x\in \X\t(\rL)\sq [X]_{\E},$ 
where $\t=\ang{\al,\ang{t_n:n\in\om}}\in\cont\al\cap\rL,$ 
$\al<\om_1.$ Then $y\not\in\X\t(\rL).$ On the other hand, 
$\t=\t(\xi)$ for some $\xi<\om_1,$ so we have $\neg\;(i)$.

Let $\ga=\al^{++}$ in $\rL,$ so that $\ga<\om_1=\Om$ and 
$\dF\al(\rL)=\dF\al(\rL_\ga).$ Then the ``term'' 
$\t'=\ang{\ga,\ang{t_n:n\in\om}}$ belongs to $\cont\ga\cap\rL,$ 
and $\X\t(\rL)=\X{\t'}(\rL_\xi)$ whenever $\ga\<\xi<\om_1.$ 
Finally, $\t'=\t(\xi)$ for some 
$\xi,$ $\ga\<\xi<\om_1,$ and then $\X\t(\rL)=\X{\t(\xi)}(\rL_\xi)$.
\qed

\subsubsection{The $\od$ subcase}
\label{closod}

We have to prove that $\E=\oE$ admits an $\od$ enumeration of the 
equivalence classes by elements of $2^{<\om_1}$. 

For every $x\in\cD,$ we define 
$\Xi(x)=\ans{\xi<\om_1:x\in [\X{\t(\xi)}(\rL)]_{\E}}$ and let 
$\phi_x\in 2^{\om_1}$ be the characteristic function of $\Xi(x).$ 
Lemma~\ref{har} implies that the $\od$ map 
$x\map\phi_x$ enumerates the \dd\E classes by elements of  
$2^{\om_1}.$ To get an enumeration by elements of 
$2^{<\om_1},$ we prove 

\ble
\label{rest}
Assume\/ $\sma\Om.$ If\/ $h\in 2^{\om_1}$ is\/ $\rod$ then 
there exists\/ $\ga<\om_1$ such that\/ $h\in\rL[h\res\ga]$.
\ele
\proof{}of the lemma. By $\sma\Om,$ there exists $\al<\om_1$ such 
that $h\in\rL[f]$ for a \dd{\col\al}generic 
over $\rL$ function $f\in \al^\om.$ Let $\uh$ be a name for 
$h$ in this forcing. 

{\em We argue in $\rL$.}
We define $H_\xi=\ans{s\in\col\al:s\,\hbox{ forces }\uh(\xi)=1}$
for all $\xi<\Om.$ (We recall that $\Om=\om_1$ in the universe but 
$\Om$ is inaccessible in $\rL$ under the assumption 
$\sma\Om$.) Since $\al<\Om,$ we have \dd{(<\Om)}many different 
sets $H_\xi.$ Therefore there exist 
an ordinal 
$\ga<\Om$ and a function $\tau:\Om\;\lra\;\ga$ such that 
$H_\xi=H_{\tau(\xi)}$ for all $\xi<\Om$. 

In the universe, this implies $h\in\rL[h\res\la],$ 
as required.\vspace{4mm}\qed

To continue the proof of the theorem, we let $\la_x$ denote the 
least ordinal $\la<\Om=\om_1$ such that 
$\rL[\phi_x]=\rL[\phi_x\res\la_x],$ for each $x\in\cD$. 

Unfortunately, the map $x\map\psi_x=\phi_x\res\la_x$ 
does not enumerate \dd\oE classes by elements of $2^{<\om_1}.$ (The 
equality $\psi_x=\psi_y$ is not sufficient for $x\E y$.) 

We utilize a more tricky idea.

Let $x\in\cD.$ Then $\psi_x=\phi_x\res\la_x\in 2^{\la_x}.$   
The set $[x]_{\oE}=\ans{x':\phi_x=\phi_{x'}}$ is $\od[\phi_x],$ 
therefore $\od[\psi_x]$ because $\phi_x\in\rL[\psi_x]$. It 
follows, by Proposition~\ref{solMb}, that $[x]_{\oE}$ includes a 
nonempty subset of the form $\X\t(\rL[\psi_x]),$ where 
$\t\in\cont{}\cap\rL[\psi_x]$.

Let $\t_x$ be the least, in the principal $\od[\psi_x]$ 
wellordering of $\rL[\psi_x],$ among the ``codes'' 
$\t\in\cont{}\cap\rL[\psi_x]$ such that 
$\emps\not=\X\t(\rL[\psi_x])\sq [x]_{\oE}$.

The map $x\map\ang{\psi_x,\t_x}$ is $\od,$ of course. 
Since the definition is \dd\oE invariant and $\E=\oE,$ we have 
$\psi_x=\psi_y$ and $\t_x=\t_y$ whenever ${x\E y}$. 

Assume now that $\psi_x=\psi_y$ and $\t_x=\t_y.$ 
In this case one and the same nonempty 
set $\X{\t_x}(\rL[\psi_x])=\X{\t_y}(\rL([\psi_y])$ is a subset 
of both $[x]_{\E}$ and $[y]_{\E},$ so $x\E y$.

Hence the map $x\map\ang{\psi_x,\t_x}$ enumerates the 
\dd\E classes by elements of the set 
\dm
\ans{\ang{\psi,\t}:\psi\in 2^{<\om_1}\,\hbox{ and }\,
\t\in\cont{}\cap\rL[\psi_x]}
\dm
This set admits an $\od$ injection in $2^{<\om_1}.$ Therefore we 
can obtain an $\od$ enumeration of the \dd\E equivalence classes 
by elements of $2^{<\om_1}.$ This ends the proof of the principal 
assertion in item \ref{1t} of Theorem~\ref{mt}. 

\subsubsection{The $\is12$ and $\is11$ subcases}
\label{closi}

Let us consider the case when $\E$ is a $\is12$ 
(resp.\ $\is11$) 
equivalence relation in item \ref{1t} of Theorem~\ref{mt}. We have 
to engineer a $\idh2$ (resp.\ $\idh1$) 
enumeration of the \dd\E equivalence classes by elements of 
$2^{<\om_1}$. 

The most natural plan would be to prove that the $\od$ enumeration 
$x\map\ang{\psi_x,\t_x}$ defined above is e.\ g. $\idh2$ provided 
$\E$ is $\is12.$ However there is no idea how to convert the 
definition of $\psi_x$ to $\idh2,$ or even to formalize it in 
$\hc.$ Fortunately we do not neet in fact the 
minimality of $\psi_x=\phi_x\res\la_x;$ all that we exploited is 
the existence of a term $\t\in\cont{}\cap\rL[\psi_x]$ such that 
$\emps\not=\X\t(\rL[\psi_x])\sq [x]_{\oE}$.

We could now define $\psi_x=\phi_x\res\la,$ where 
$\la=\la_x$ is the least ordinal $\la<\om_1$ such that 
$\cont{}\cap\rL[\psi_x]$ contains the required term. This can be 
formalized in $\hc,$ but hardly as a $\idh2$ definition: indeed, 
e.\ g. the condition $\X\t(\rL[\psi_x])\sq [x]_{\oE}$ does not 
look like better than $\iph2$.

The correct plan includes one more idea, originally due to Hjorth 
and Kechris~\cite{hk}: certain requirements 
are eliminated by forcing over a countable submodel. 

Let us consider details. We recall that $\sma\Om$ is assumed. 

Let $x\in\cD.$ We define $\vpi_x\in 2^{\om_1}$ 
by $\vpi_x(\xi)=1$ iff $x\in[\X{\tk\xi}(\rL_\xi)]_{\E}$ for all 
$\xi<\om_1.$ (Pay attention on the similarity and the difference 
between $\vpi_x$ and $\phi_x$ above.) 

\bdf
\label{psi}
We let $T_{\E}$ be the set of all triples $\ang{x,\psi,t}$ such that 
$x\in\cD,$ $\psi\in 2^{<\om_1},$  
$\t\in\cont\al\cap\rL_{\ga}[\psi],$ where 
$\al<\ga=\dom\psi<\om_1,$ and the following conditions \ref{aa} 
through \ref{cc} are satisfied. 
\ben
\def\theenumi{\hskip2pt\rmt(\alph{enumi})\hskip2pt}
\def\labelenumi{\theenumi} 
\itla{aa} \msur $\rL_{\ga}[\psi]$ models $\ZFC^-$ (minus power 
set) so that $\psi$ can occur as an extra class parameter in 
Replacement and Separation.\its

\itla{bb}
It is true in $\rL_{\ga}[\psi]$ that $\ang{\La,\La}$ 
forces ${\C\t(\unf)\E\C\t(\ung})$ in the sense of 
$\col\al\!\ti\!\col\al$ as the forcing, where $\unf$ and 
$\ung$ are the names for the generic functions in $\al^\om$.\its

\itla{xx}\msur 
$\psi=\vpi_x\res\ga$.\its

\itla{cc} 
\msur 
$x$ belongs to $[\X\t(\rL_\ga[\psi])]_{\E}$. \its
\een
A point $x\in\cD$ is {\em \dd\E classifiable} iff there 
exist $\psi$ and $\t$ such that $\ang{x,\psi,\t}\in T_{\E}$.\qed
\edf

\ble
\label{key2}
Assume\/ $\sma\Om.$ If\/ $\E$ is a\/ $\is12$ equivalence and\/ 
$\E=\oE$ then all points\/ $x\in\cD$ are\/ \dd\E classifiable.
\ele
\proof Let $x\in\cD$. Then $\vpi_x$ is $\od[x],$ 
therefore $\vpi_x\in\rL[x]$ by Proposition~\ref{solM}. 
Since $\E=\oE,$ Lemma~\ref{har} implies that the set $[x]_{\E}$ 
is $\od[\vpi_x].$ Therefore by Proposition~\ref{solMb} we have 
$x\in\X\t(\rL[\vpi_x])\sq [x]_{\E}$ for some 
$\t\in\cont\al\cap\rL[\vpi_x],$ $\al<\Om=\om_1$.

The model $\rL_{\om_1}[\vpi_x]$ has an elementary submodel 
$\rL_\ga[\psi],$ where $\ga<\om_1$ and $\psi=\vpi_x\res\ga,$ 
containing $\t$ and $\al.$ We prove that 
$\ang{x,\psi,\t}\in T_{\E}.$ Since conditions \ref{aa} and 
\ref{xx} of Definition~\ref{psi} obviously hold 
for $\rL_\ga[\psi],$ let us check requirements \ref{bb}, 
\ref{cc}.\vspace{1mm}

{\em We check\/ \ref{bb}.} 
Indeed otherwise there exist conditions $u,\,v\in\col\al$ such 
that $\ang{u,v}$ forces ${\C{\t}(\unf)\nE \C{\t}(\ung)}$ in 
$\rL_{\ga}[\psi]$ in the sense of $\col\al\!\ti\!\col\al$ as the 
notion of forcing. Then $\ang{u,v}$ also forces 
${\C{\t}(\unf)\nE \C{\t}(\ung)}$ in $\rL_{\om_1}[\vpi_x]$. 
Let us consider an 
\dd{\col\al\!\ti\!\col\al}generic over $\rL[\vpi_x]$ pair 
$\ang{f,g}\in \al^\om\ti\al^\om$ such that $u\subset f$ and 
$v\subset g.$ Then both $y=\C\t(f)$ and $z=\C\t(g)$ belong to 
$\X \t(\rL[\vpi_x]),$ so ${y\E z}$ because 
$\X\t(\rL[\vpi_x])\sq [x]_{\E}$.

Notice that $\ang{f,g}$ also is generic over 
$\rL_{\om_1}[\vpi_x].$ We observe that ${y\E z}$ is 
{\em false\/} in $\rL_{\om_1}[\vpi_x,f,g],$ that is, in 
$\rL[\vpi_x,f,g],$ by the choice of $u$ and $v.$ But ${y\E z}$ is 
a $\is12$ formula, therefore absolute for transitive models 
containing all ordinals, contradiction.\vspace{1mm}


{\em We check\/ \ref{cc}.} Take any \dd{\col\al}generic 
over $\rL[\vpi_x]$ function $f\in\al^\om.$ Then $y=\C\t(f)$ belongs 
to $\X\t(\rL[\vpi_x]),$ hence ${y\E x}.$ On the other hand, $f$ is 
generic over $\rL_{\ga}[\psi]$.\vspace{1mm}

Thus $\ang{x,\psi,\t}\in T_{\E}.$ This means that 
$x$ is \dd\E classifiable, as required.\qed\vspace{4mm}

Thus, for each \dd\E classifiable $x,$ all countable sequences 
$\psi$ which satisfy $T_{\E}(x,\psi,\t)$ for some $\t,$ are restrictions 
of one and the same sequence $\vpi_x\in 2^{\om_1},$ 
defined 
above.

\bdf
\label{U}
Let $x\in\cD.$ It follows from Lemma~\ref{key2} that there 
exists the least ordinal $\ga=\ga_x<\om_1$ 
such that $T_{\E}(x,\vpi_x\res\ga,\t)$ for some $\t.$ We put 
$\psi_x=\vpi_x\res\ga$ and let $\t_x$ denote the least, in 
the sense of the $\od[\psi_x]$ wellordering of $\rL_{\ga}[\psi_x],$ 
``term'' 
$\t\in\cont{}[\psi_x]\cap \rL_{\ga}[\psi_x]$ which satisfies 
$T_{\E}(x,\psi_x,\t).$ We put $U(x)=\ang{\psi_x,\t_x}$.\qed 
\edf

\ble
\label{inv}
Assume\/ $\sma\Om.$ If\/ $\E$ is a\/ $\is12$ equivalence and\/ 
$\E=\oE$ then the map\/ $U$ enumerates the\/ \dd\E classes. 
\ele
\proof If $x\E y$ then $U(x)=U(y)$ because Definition~\ref{psi} 
is \dd\E invariant for $x$. 

Let us prove the converse. Assume that $U(x)=U(y),$ that is, in 
particular, $\psi_x=\psi_y=\psi\in 2^{<\om}$ and 
${\t_x=\t_y=\t\in\cont\al[\psi]\cap\rL_{\ga}[\psi],}$ 
where $\al<\ga=\dom\psi<\om_1=\Om$. 



By \ref{cc} we have ${\C\t(f)\E x}$ and ${\C\t(g)\E y}$ for 
some \dd{\col\al}generic over $\rL_{\ga}[\psi]$ functions 
$f,\,g\in\al^\om.$ Let us consider an \dd{\col\al}generic over 
both $\rL_{\ga}[\psi,f]$ and $\rL_{\ga}[\psi,g]$ function 
$h\in\al^\om.$ Then, by \ref{bb}, ${\C\t(h)\E\C\t(f)}$ holds in 
$\rL_{\ga}[\psi,f,h],$ therefore in the universe because $\E$ is 
$\is12.$ Similarly, we have ${\C\t(h)\E\C\t(g)}.$ It follows that 
${\C\t(f)\E\C\t(g)},$ hence ${x\E y},$ as required.
%
\qed

\ble
\label{def}
Suppose that\/ $\E$ is\/ $\is12$ {\rmt({\it resp. $\is11$})} 
and\/ $\E=\oE.$ 
Then\/ $U$ is a function of class\/ $\idh2$ 
{\rmt({\it resp. $\idh1$}).} 
\ele
\proof It suffices to check that the set $T_{\E}$ is $\idh2$ (resp. 
$\idh1$). 

Notice that conditions \ref{aa} and \ref{bb} in 
Definition~\ref{psi} are $\idh1$ because they reflect truth 
within $\rL_\ga[\psi]$ and the enumeration $\tk\xi$ was chosen 
in $\idh1$. 

Suppose that $\E$ is $\is12,$ that is, $\ish1.$ Then condition 
\ref{cc} is obviously $\ish1.$ Condition \ref{xx} can be converted 
to $\idh2$ (in fact a bounded quantifier $\forall\,\ba<\ga$ over a 
conjunction of $\ish1$ and $\iph1$ relations). Indeed we observe 
that \ref{xx} is equivalent to 
\dm
\forall\,\xi<\ga\;(\psi(\xi)=1\;\,\llra\;\,
x\in [\X{\t(\xi)}(\rL_\xi)]_{\E}) \eqno{(\ast)}
\dm
by Lemma~\ref{har}. 

The case when $\E$ belongs to $\is11$ is more difficult. 

Let us first consider condition~\ref{cc}. Immediately, it is 
$\is12,$ therefore $\ish1,$ so it remains to 
convert it also to a $\iph1$ form. Notice that in the assumption 
of \ref{aa} and \ref{bb}, the set $X=\X\t(\rL_{\ga}[\psi])$ 
consists of pairwise \dd\E equivalent points: this was actually 
showed in the proof of Lemma~\ref{inv}. Therefore, since 
obviously $\X\t(\rL_{\ga}[\psi])\not=\emps,$ \ref{cc} is 
equivalent to 
$\forall\,y\in \X\t\;(\rL_{\ga}[\psi])\;(x\E y).$ This is 
clearly $\iph1$ provided $\E$ is $\ip12$. 

Let us consider \ref{xx}. The right--hand side of the equivalence 
$(\ast)$ is $\is11$ with inserted $\idh1$ functions, therefore 
$\idh1.$ It follows that $(\ast)$ itself is $\idh1,$ as required.
\qed\vspace{4mm}

This completes the proof of the additional part ($\is11$ 
and $\is12$ relations) in item~\ref{1t} of Theorem~\ref{mt}.

\newpage
\subsection{$\od$ topology and the forcing}
\label{top}

This section starts the proof of item \ref{2t} of Theorem 
\ref{mt} for a given $\od$ equivalence relation $\E$ in the 
assumption $\sma\Om$.

We have to embed $\Eo$ in $\E$ continuously. The embedding will 
be defined in the next section; here we obtain some useful 
preliminary results related to the defined above topology $\top,$ 
an associated forcing and the relevant product forcing. At the 
end of the section, we introduce the set $H$ of all points 
$x\in\cD$ whose \dd\oE classes are bigger than \dd\E classes; 
$H$ is nonempty in the assumption $\E\sneq\oE$.

The reasoning is based on special properties of the $\od$ topology 
$\top,$ having a semblance of the Gandy -- Harrington topology 
(even in a simplified form because some specific $\is11$ details 
vanish). In particular, the topology is strongly Choquet. However 
we shall not utilize this property (and shall not prove it). We 
take indeed another way. The reasoning will be organized as a 
sequence of straight forcing arguments. This manner of treatment 
of equivalence relations was taken from Miller~\cite{mill}.

\subsubsection{Topology and the forcing}
\label{tt}
 
The topology $\top$ obviously 
does not have a countable base; but it has one in 
a local sense. A set $X$ will be called \dd\top{\it separable\/} 
if the $\od$ power set ${\oP(X)=\cP(X)\cap\od}$ has only 
countably many different $\od$ subsets. 

\ble
\label{dizl}
Assume\/ $\sma\Om.$ Let\/ $\al<\Om$ and\/ $\t\in\cont\al\cap\rL.$ 
Then\/ $X=\X\t(\rL)$ is\/ \dd\top separable.
\ele
\proof By Proposition~\ref{solMb} every $\od$ subset of $X$ is 
uniquely determined by an $\od$ subset of $\col\al.$ Since each 
$\od$ set $S\sq\col\al$ is constructible (Proposition~\ref{solM}), 
we obtain an $\od$ map $h:\al^+\,\hbox{ onto }\,\oP(X),$ 
where $\al^+$ is the least cardinal in $\rL$ bigger than $\al.$ 
Therefore $\oP(X)$ has \dd{\<\al^{++}}many $\od$ subsets. It 
remains to notice that $\al^{++}<\Om$ because $\Om$ is 
inaccessible in $\rL$.\qed\vspace{4mm}

Let $\dX=\ans{X\sq\cD:X\,\hbox{ is }\,\od\;\hbox{ and nonempty}\,}$.

Let us consider $\dX$ as a forcing notion (smaller sets are 
stronger conditions) for generic extensions of $\rL$ in the 
assumption $\sma\Om.$ Of course formally $\dX\not\in\rL,$ but 
$\dX$ is $\od$ order isomorphic to 
a partially ordered set in $\rL$. (Indeed it is known that there 
exists an $\od$ map $\phi:$ ordinals onto the class of all $\od$ 
sets. Since $\dX$ itself is $\od,$ $\dX$ is a 1--1 image of an 
$\od$ set $\dX'$ of ordinals via $\phi.$ By Proposition~\ref{solM} 
both $\dX'$ and the \dd{\phi\hspace{0.5pt}}preimage of the order 
on $\dX$ belong to $\rL$.)

It also is true that a set $G\sq\dX$ is \dd\dX generic over $\rL$ 
iff it nonempty intersects every dense $\od$ subset of $\dX$. 

\bcor
\label{exis}
Assume\/ $\sma\Om.$ If a set\/ ${X\in\dX}$ is nonempty then 
there exists an\/ \dd\dX generic over\/ $\rL$ set\/ $G\sq\dX$ 
containing\/ $X$.
\ecor
\proof We can suppose, by Proposition~\ref{solMb}, that 
$X=\X\t(\rL)$ where $\t\in\cont\al\cap\rL$ and $\al<\Om.$ 
Now apply Lemma~\ref{dizl}.\qed

\ble 
\label{choq-cor}
Assume\/ $\sma\Om.$ If\/ $G\sq\dX$ is a generic over\/ $\rL$ 
set then the intersection\/ $\bigcap G$ is a singleton 
$\ans{a}=\ans{a_G}$.
\ele
\proof Assume that this is not the case. Let $\dX'\in\rL$ be a 
constructible p. o. set order isomorphic $\dX$ via an $\od$ 
function $f:\dX'\,\hbox{ onto }\,\dX.$ Then $G'=\phi^{-1}(G)$ is 
\dd{\dX'}generic over $\rL.$ We assert that the statement that 
$\bigcap G$ is not a singleton can be converted to a sentence 
relativized to $\rL[G']$. 

(Indeed, it follows from the reasoning in the proof of 
Lemma~\ref{dizl} that $\rL[G']$ is in fact a \dd Pgeneric 
extension of $\rL$ for a certain set $P\in\rL,$ $P\sq\dX'$ of 
a cardinality $\al<\Om$ in $\rL.$ The next \dd\rL cardinal 
$\al^+$ is $<\Om.$ Therefore $G'$ belongs to a 
\dd{\col{\al^+}}generic extension of $\rL,$ so $G'$ is \dd\Om weak 
over $\rL.$ Then by Proposition \ref{solM} the universe 
is a Solovay extension of $\rL[G'].$ 
This is enough to convert any statement about $G'$ in $\rV$ -- 
like the statement: $\bigcap \phi\ima G'$ is not a singleton -- 
to a sentence relativized to $\rL[G']$.)

Then there exists ${X\in \dX}$ such that $\bigcap G$ is 
{\em not\/} a singleton for {\em every\/} generic over $\rL$ set 
$G\sq \dX$ containing $X.$ We can assume that 
$X=\X\t(\rL),$ where ${\t\in\cont\al\cap\rL},$ $\al<\Om.$ Then $X$ 
is \dd\top separable; let $\ans{\cX_n:n\in\om}$ be an enumeration 
of all $\od$ dense subsets of $\oP(X).$ Using 
Proposition~\ref{solMb} (item~\ref{sm6}), we obtain an increasing 
\dd{\col\al}generic over $\rL$ sequence 
$u_0\sq u_1\sq u_2\sq...$ of $u_n\in\col\al$ such that 
$X_n=\X{\t\hspace{0.1em}{u_n}}(\rL)\in\cX_n.$ Obviously this 
gives an \dd\dX generic over $\rL$ set $G\sq\dX$ containing $X$ 
and all $X_n$.

Now let $f=\bigcup_{n\in\om}u_n;$ $f\in\al^\om$ and $f$ is 
\dd{\col\al}generic over $\rL.$ Then $x=\C\t(f)\in X_n$ for all 
$n,$ so $x\in\bigcap G.$ Since $\bigcap G$ obviously cannot 
contain more than one point, it is a singleton, so we get a 
contradiction with the choice of $X$.\qed\vspace{4mm}

Reals $a_G$ will be called \dd\od{\it generic over\/} $\rL$.

\subsubsection{The product forcing}

We recall that $\E$ is an $\od$ equivalence on 
$\cD$ and $\oE$ is the \dd{\top^2}closure of $\E$.

For a set ${P\sq\cD^2,}$ we put 
${\pri P=\ans{x:\exists\,y\;P(x,y)}}$ and 
${\prii P=\ans{y:\exists\,x\;P(x,y)}.}$ Notice that if $P$ is 
$\od,$ so are $\pri P$ and $\prii P$. 

The classical reasoning in Harrington, Kechris, and 
Louveau~\cite{hkl} plays on interactions between $\E$ and $\oE.$ 
In the forcing setting, we have to fix a restriction by $\oE$ 
directly in the definition of the product forcing. Thus we 
consider 
\dm
\dP=\dP(\oE)=\ans{P\sq\oE:P\,\hbox{ is }\od\;\hbox{ and 
nonempty and }\,P=(\pri P\ti \prii P)\cap\oE}
\dm 
as a forcing notion. As above for $\dX,$ the fact that formally 
$\dP$ does not belong to $\rL$ does not cause essential problems.

The following assertion connects $\dP$ and $\dX$.

\bass
\label{proe}
Assume $\sma\Om.$ Then\its
\ben
\def\theenumi{\rmt{\arabic{enumi}}}
\item If\/ $P\in\dP$ then\/ $\pri P$ and\/ $\prii P$ belong to 
$\dX$.\its

\item If\/ $X,\,Y\in\dX$ and\/ $P=(X\ti Y)\cap \oE\not=\emps$ 
then $P\in\dP$.\its

\itla{i3} 
If\/ $P\in\dP,$ $X\in\dX,$ $X\sq\pri P,$ then 
there exists\/ $Q\in\dP,$ $Q\sq P,$ such that $X=\pri Q.$  
Similarly for $\prii$.
\een
\eass
\proof Set $Q=\ans{\ang{x,y}\in P:x\in X\cj y\oE x}$ in 
item~\ref{i3}.\qed\vspace{4mm}

A set $P\in\dP$ is {\em \dd\dP separable} if 
the set $\dP_{\sq P}=\ans{Q\in\dP:Q\sq P}$ has only countably many 
different $\od$ subsets. 

\ble 
\label{dizl2}
Assume\/ $\sma\Om.$ Suppose that\/ ${P=(X\ti Y)\cap\oE}$ 
is nonempty, where\/ $X=\X\t(\rL),$ $Y=\X{\t'}(\rL),$ and\/ 
$\t,\,\t'\in\cont{}\cap\rL.$ 
Then\/ $P\in\dP$ and\/ $P$ is\/ \dd\dP separable.
\ele
\proof We have $P\in\dP$ by Assertion~\ref{proe}. A proof of the 
\dd\dP separability can be obtained by a minor modification of 
the proof of Lemma~\ref{dizl}.\qed

\ble
\label{dp2oe}
Assume\/ $\sma\Om.$ Let\/ $G\sq \dP$ be a\/ \dd{\dP}generic over\/ 
$\rL$ set. Then the intersection\/ $\bigcap G$ contains a single 
point\/ $\ang{a,b}$ where\/ $a$ and\/ $b$ are\/ \dd\od generic 
over $\rL$ and $a\oE b$.
\ele
\proof By Assertion~\ref{proe}, both $G_1=\ans{\pri P:P\in G}$ and 
$G_2=\ans{\pri P:P\in G}$ are \dd\od generic over $\rL$ subsets of 
$\dX,$ so that there exist unique \dd\od generic over $\rL$ points 
$a=a_{G_1}$ and $b=a_{G_2}.$ It remains to show that 
$\ang{a,b}\in\oE$.

Suppose not. There exists an \dd\E invariant $\od$ set 
$A$ such that we have $x\in A$ and $y\in B=\cD\setminus A.$ Then 
$A\in G_1$ and $B\in G_2$ by the genericity. There exists a 
condition $P\in G$ such that $\pri P\sq A$ and $\prii B\sq B,$ 
therefore ${P\sq (A\ti B)\cap\oE=\emps},$ which is 
impossible.\qed\vspace{4mm}

Pairs $\ang{a,b}$ as in Lemma~\ref{dp2oe} will be called 
\dd\dP{\it generic\/} and denoted by $\ang{a_G,b_G}$.

For sets $X$ and $Y$ and a binary relation $\R{}\,,$ let us write 
${X\R{}Y}$ if and only if 
$\forall\,x\in X\;\exists\,y\in Y\;(x\R{} y)$ \ and \ 
$\forall\,y\in Y\;\exists\,x\in X\;(x\R{} y)$.

\ble
\label{1for2}
Assume\/ $\sma\Om.$ Suppose that\/ $P_0\in\dP,$ points\/ 
$a,\,a'\in X_0=\pri P_0$ are\/ 
\dd\od generic over\/ $\rL,$ and\/ ${a\oE a'.}$ There exists a 
point\/ $b$ such that both\/ $\ang{a,b}$ and $\ang{a',b}$ belong 
to\/ $P_0$ and are\/ \dd{\dP}generic pairs.
\ele
\proof By Lemma \ref{dizl2} and Proposition~\ref{solMb} there 
exists a \dd\dP separable set $P_1\sq P_0$ such that 
$a\in X_1=\pri P_1.$ We put $Y_1=\prii P_1;$ then $X_1\oE Y_1,$ 
and $P_1=(X_1\ti Y_1)\cap\oE$. 

We let $P'=\ans{\ang{x,y}\in P_0:y\in Y_1}.$ Then $P'\in \dP$ 
and $P_1\sq P'\sq P_0.$ Furthermore $a'\in X'=\pri P'.$ (Indeed, 
since ${a\in X_1}$ and ${X_1\oE Y_1},$ there exists $y\in Y_1$ 
such that $a\oE y;$ then $a'\oE y$ as well because $a\oE a',$ 
therefore $\ang{a',y}\in P'$.) As above  
there exists a \dd\dP separable set $P'_1\sq P'$ such that 
$a'\in X'_1=\pri P'_1.$ Then $Y'_1=\prii P'_1\sq Y_1$. 

It follows from the choice of $P$ and $P'$ that $\dP$ admits only 
countably many different dense $\od$ sets below $P_1$ and below 
$P'_1.$ Let $\ans{\cP_n:n\in\om}$ and $\ans{\cP'_n:n\in\om}$ 
be\pagebreak[3] enumerations of both families of dense 
sets. We define sets $P_n,\,P'_n\in\dP\;\;(n\in\om)$ 
satisfying the following conditions:\its
\ben
\def\theenumi{(\roman{enumi})}
\def\labelenumi{\theenumi}
\itla{i}
$a\in X_n=\pri P_n$ \ and \ $a'\in X'_n=\pri P'_n$;\its

\itla{ii}
$Y'_n=\prii P'_n \sq Y_n=\prii P_n$ \ and \ $Y_{n+1}\sq Y'_n$;\its

\itla{iii}
$P_{n+1}\sq P_n\,,\,$ $P'_{n+1}\sq P'_n\,,\,$ 
$P_n\in \cP_{n-2}\,,\,$ and \ $P'_n\in \cP'_{n-2}$.\its
\een
By \ref{iii} both sequences $\ans{P_n:n\in\om}$ and 
$\ans{P'_n:n\in\om}$ are \dd\dP generic over $\rL,$ so by 
Lemma~\ref{dp2oe} they result in two generic pairs,  
$\ang{a,b}\in P_0$ and $\ang{a',b}\in P_0, $ having the first 
terms equal to $a$ and $a'$ by \ref{i} and second terms equal to 
each other by \ref{ii}. Thus 
it suffices to conduct the construction of $P_n$ and $P'_n$.

The construction goes on by induction on $n$.

Assume that $P_n$ and $P'_n$ have been defined. We define 
$P_{n+1}.$ By~\ref{ii} and Assertion~\ref{proe}, the set 
${P=(X_n\ti Y'_n)\cap\oE\sq P_n}$ belongs to $\dP$ and 
$a\in X=\pri P.$ (Indeed, $\ang{a,b}\in P,$ where $b$ satisfies 
$\ang{a',b}\in P'_n,$ because ${a\oE a'}$.) However $\cP_{n-1}$ 
is dense in $\dP$ below $P\sq P_0;$ therefore 
${\pri \cP_{n-1}=\ans{\pri P':P'\in \cP_{n-1}}}$ is dense in $\dX$ 
below\pagebreak[3] 
$X=\pri P.$ Since $a$ is generic, we have $a\in \pri P'$ for 
some $P'\in \cP_{n-1},$ $P'\sq P.$ It\pagebreak[3] 
remains to put $P_{n+1}=P',$ 
and then $X_{n+1}=\pri P_{n+1}$ and $Y_{n+1}=\prii P_{n+1}$.

After this, to define $P'_{n+1}$ we let 
$P=(X'_n\ti Y_{n+1})\cap\oE,$ etc.\qed 

\subsubsection{The key set}
\label{second}

We recall that $\sma\Om$ is assumed, $\E$ is an $\od$ 
equivalence on $\cD,$ and $\oE$ is the \hbox{\dd{\top^2}closure} 
of $\E$ in $\cD^2.$ We shall also suppose that $\E\sneq\oE,$ in 
accordance to item~\ref{2t} of Theorem~\ref{mt}. Then there exist 
\dd\oE classes which include more than one \dd\E class. (In fact 
we shall have no other use of the hypothesis $\E\sneq\oE$.) 
We define the union of all those \dd\oE classes, 
\dm
H=\ans{x\in\cD:\exists\,y\in\cD\;(x\oE y\cj x\nE y)}\,,
\dm
the ``key set'' from the title. The role of this set in the 
reasoning below is entirely similar to the role of the 
corresponding set $V$ in Harrington, Kechris, and 
Louveau~\cite{hkl}. 

Obviously $H$ is $\od,$ nonempty, and \dd\E invariant. 
Furthermore, $H'=H^2\cap\oE\not=\emps$ (in fact both projections 
of $\cH$ are equal to $H$),  
so that in particular $H'\in\dP$ by Assertion~\ref{proe}. 

\ble
\label{noE} 
Assume\/ $\sma\Om.$ If\/ $a,b\in H$ and\/ $\ang{a,b}$ is\/ 
\dd\dP generic over $\rL$ then ${a\nE b}\hspace{1.5pt}.$
\ele
\proof Otherwise there exists a set $P\in\dP,$ $P\sq H\ti H$ such 
that $a\E b$ holds for {\it all\/} \dd\dP generic $\ang{a,b}\in P.$ 
We conclude that then $a\oE a'\;\lra\;a\E a'$ for all \dd\od generic 
points $a,\,a'\in X=\pri P;$ indeed, take $b$ such that both 
$\ang{a,b}\in P$ and $\ang{a',b}\in P$ are \dd\dP generic, 
by Lemma~\ref{1for2}. In other words the relations $\E$ and 
$\oE$ coincide on the set 
${Y=\ans{x\in X:x\,\hbox{ is \dd\od generic over }\,\rL}\in\dX.}$ 
($Y\not=\emps$ by corollaries \ref{exis} and \ref{choq-cor}.) 

Moreover, $\E$ and $\oE$ coincide on the set 
$Z=[Y]_{\E}.$ Indeed if $z,\,z'\in Z,$ ${z\oE z'},$ 
then let ${y,\,y'\in Y}$ satisfy ${z\E y}$ and ${z'\E y'}.$ 
Then ${y\oE y'},$ therefore ${y\E y'},$ which implies $z\E z'.$  

We conclude that $Y\cap H=\emps$. 

(Indeed, suppose that 
$x\in Y\cap H.$ Then by definition there exists $y\in\cD$ 
such that ${x\oE y}$ but ${x\nE y}.$ Then ${y\not\in Z}$ because 
$\E$ and $\oE$ coincide on $Z.$ Thus the pair $\ang{x,y}$ belongs 
to the $\od$ set $P=Y\ti (\cD\setminus Z).$ Notice that $P$ 
does not intersect $\E$ by definition of $Z.$ Therefore 
$\ang{x,y}$ cannot belong to the closure $\oE$ of $\E,$ 
contradiction.) 

But $\emps\not=Y\sq X\sq H,$ contradiction.\qed\vspace{4mm}

Lemma~\ref{noE} is a counterpart of the proposition in 
Harrington, Kechris, Louveau~\cite{hkl} that $\E\res H$ is 
meager in $\oE\res H.$ But in fact the main content of this 
argument in~\cite{hkl} was implicitly taken by Lemma~\ref{1for2}. 

\ble
\label{E}
Assume\/ $\sma\Om.$ Let\/ $X,\,Y\sq H$ be nonempty\/ $\od$ sets 
and\/ ${X\oE Y}.$ There 
exist nonempty\/ $\od$ sets\/ $X'\sq X$ and\/ $Y'\sq Y$ such 
that\/ $X'\cap Y'=\emps$ but still\/ $X'\oE Y'$.
\ele
\proof There exist points $x_0\in X$ and $y_0\in Y$ such that 
$x_0\not= y_0$ but ${x_0\oE y_0}.$ (Otherwise $X=Y,$ and $\oE$ is 
the equality on $X,$ which is impossible, see the previous proof.) 
Let $U$ and $V$ be disjoint Baire intervals in $\cD$ containing 
resp. $x_0$ and $y_0.$ The sets $X'= X\cap U \cap [Y\cap V]_{\oE}$ 
and $Y'= Y\cap V \cap [X\cap U]_{\oE}$ are as required.\qed

\newpage

\subsection{The embedding}
\label{or}

In this section we accomplish the proof of Theorem~\ref{mt}, 
therefore Theorem~\ref{main} (see Section~\ref{appr}). Thus we 
prove, assuming $\sma\Om$ and $\E\sneq\oE,$ that $\E,$ the given 
$\od$ equivalence on $\cD,$ continuously embeds $\Eo$. 

\subsubsection{The embedding}
\label{embed}

By the assumption the set $H$ of Subsection~\ref{second} is 
nonempty; obviously $H$ is $\od.$ By lemmas \ref{dizl} and 
\ref{dizl2} there exists a nonempty \dd\top separable $\od$ set 
$X_0\sq H$ such that the set ${P_0=(X_0\ti X_0)\cap\oE}$ belongs 
to $\dP$ and is \dd\dP separable; 
$\pri P_0=\prii P_0=X_0$.  

We define a family of sets $X_u\;\;(u\in 2^{<\om})$ satisfying 
\its
\ben
\def\theenumi{(\alph{enumi})}
\def\labelenumi{\theenumi}
\itla{a} 
$X_u\sq X_0,$ $X_u$ is nonempty and $\od,$ and $X_{u\we i}\sq X_u,$ 
for all $u$ and $i$.\its
\een
In addition to the sets $X_u,$ we shall define sets 
$\J{uv}$ for {\em some} pairs $\ang{u,v},$ to provide important 
interconnections between branches in $2^{<\om}.$ 

Let $u,\,v\in 2^n.$ We say that $\ang{u,v}$ is a {\em neighbouring 
pair\/} iff $u=0^k\we 0\we r$ and $v=0^k\we 1\we r$ for some $k<n$ 
($0^k$ is the sequence of $k$ terms each equal to $0$) and some 
$r\in 2^{n-k-1}$ (possibly $k=n-1,$ that is, $r=\La$). 

Thus we define sets $\J{uv}\sq X_u\ti X_v$ for all neighbouring pairs 
$\ang{u,v},$ so that the following requirements \ref{b} and  
\ref{d} will be satisfied.\its
\ben
\def\theenumi{(\alph{enumi})}
\def\labelenumi{\theenumi}
\setcounter{enumi}{1}
\itla{b} \msur 
$\J{uv}$ is $\od,$ $\pri \J{uv}=X_u,$ $\prii \J{uv}=X_v,$ and 
$\J{u\we i\,,\,v\we i}\sq \J{uv}$ for every neighbouring pair 
$\ang{u,v}$ and each $i\in\ans{0,1}$.\its

\itla{d} 
For any $k,$ the set $\J k=\J{0^k\we 0\,,\,0^k\we 1}$ is 
\dd\top separable, and $\J k\sq \E$.\its
\een
Notice that if $\ang{u,v}$ is neighbouring then 
$\ang{u\we i,v\we i}$ is neighbouring, but $\ang{u\we i,v\we j}$ 
is not neighbouring for $i\not=j$ (unless $u=v=0^k$ for some $k$). 

This implies $X_u \J{uv} X_v,$ therefore 
$X_u\E X_v,$ for all neighbouring pairs $u,\,v.$~\footnote
{\ We recall that $X\J{}Y$ means that 
$\forall\,x\in X\;\exists\,y\in Y\;(x\J{} y)$ and 
$\forall\,y\in Y\;\exists\,x\in X\;(x\J{} y)$.}

\brem
\label{newrem}
Every pair of $u,\,v\in 2^n$ can be tied in $2^n$ by a 
finite chain of neighbouring pairs. It follows that 
${X_u\E X_v}$ and ${X_u\oE X_v}$ hold for {\em all} pairs 
$u,\,v\in 2^n$.\qed
\erem

Three more requirements will concern genericity. 

Let $\ans{\cX_n:n\in\om}$ be a fixed (not necessarily $\od$) 
enumeration of all dense in $\dX$ below $X_0$ subsets of $\dX.$ 
Let $\ans{\cP_n:n\in\om}$ be a fixed enumeration of all dense in 
$\dP$ below $P_0$ subsets of $\dP.$ It is assumed that 
$\cX_{n+1}\sq\cX_n$ and $\cP_{n+1}\sq\cP_n.$ Note that 
${\cX'=\ans{P\in\dP: P\sq P_0\cj \pri P\cap\prii P=\emptyset}}$  
is dense in $\dP$ below $P_0$ by Lemma~\ref{E}, so we can suppose 
in addition that $\cP_0=\cX'$. 

In general, for any \dd\top separable set $S$ let 
$\ans{\cX_n(S):n\in\om}$ be a fixed enumeration of all 
dense subsets in the algebra $\oP(S)\setminus\ans{\emps}.$ 
It is assumed that $\cX_{n+1}(S)\sq\cX_n(S)$.\its
\ben
\def\theenumi{({\rmt g}\arabic{enumi})}
\def\labelenumi{\theenumi}
\itla{g1} \msur 
$X_u\in \cX_n$ whenever $u\in 2^n$.\its

\itla{g2}
If $u,\,v\in 2^n$ and $u(n\mins 1)\not=v(n\mins 1)$ (that is, the 
last terms of $u,\,v$ are different), then 
$P_{uv}=(X_u\ti X_v)\cap\oE\in \cP_n.$ \hfill --- \hfill 
{\it In fact this implies\/} \ref{g1}.\its

\itla{g3} 
If $\ang{u,v}=\ang{0^k\we 0\we r,0^k\we 1\we r}\in (2^n)^2$ 
then $\J{uv}\in \cX_n(\J k)$.\its
\een
In particular \ref{g1} implies by Corollary~\ref{choq-cor} that 
for any $a\in 2^\om$ the intersection 
$\bigcap_{n\in\om}X_{a\res n}$ contains a single point, denoted 
by $\phi(a),$ which is \dd\dX generic over $\rL,$ and the map 
$\phi$ is continuous in the Polish sense. 

\bass
\label{nass}
$\phi$ is a continuous 1--1 embedding\/ $\Eo$ 
into $\E$.
\eass
\proof Let us prove that $\phi$ is 1--1. Suppose that 
${a\not=b\in 2^\om.}$ Then ${a(n\mins 1)\not=b(n\mins 1)}$ for 
some $n.$ Let ${u=a\res n},$ ${v=b\res n},$ 
so that we have $x=\phi(a)\in X_u$ and $y=\phi(b)\in X_v.$ But 
then the set ${P=(X_u\ti X_v)\cap \oE}$ belongs to $\cP_n$ by 
\ref{g2}, therefore to $\cP_0.$ This implies 
$X_u\cap X_v=\emptyset$ by definition of $\cP_0,$ 
hence $\phi(a)\not=\phi(b)$ as required.  

Furthermore if $a\nEo b$ (which means that $a(k)\not=b(k)$ for 
infinitely many numbers $k$) then $\ang{\phi(a),\phi(b)}$ is 
\dd\dP generic by \ref{g2}, so $\phi(a)\nE \phi(b)$ by 
Lemma~\ref{noE}.

Let us finally verify that ${a\Eo b}$ implies 
${\phi(a)\E \phi(b)}.$
It is sufficient to prove that 
${\phi(0^k\we 0\we c)\E \phi(0^k\we 1\we c)}$ holds for all 
${k\in\om}$ and ${c\in 2^\om,}$ simply because every pair of 
$u,\,v \in 2^n$ is tied in $2^n$ by a chain
of neighbouring pairs. 

The sequence of sets 
$W_m=\J{0^k\we 0\we c\res m\,,\,0^k\we 1\we c\res m}$ $(m\in\om)$ 
is then generic over $\rL$ by \ref{g3} in the sense of the forcing 
$\oP(\J k)\setminus\ans{\emps}$ (we recall that 
$\J k=\J{0^k\we 0\,,\,0^k\we 1}$), 
which is simply a copy of $\dX,$ so that by 
Corollary~\ref{choq-cor} the intersection of all sets $W_m$ is a 
singleton. Obviously the singleton can be only equal to 
$\ang{\phi(0^k\we 0\we c)\,,\,\phi(0^k\we 1\we c)}.$ We conclude 
that $\phi(0^k\we 0\we c)\E \phi(0^k\we 1\we c),$ as required.\qed

\subsubsection{Restriction lemma}

Thus the theorem is reduced to the construction of sets $X_u$ 
and $\J{uv}$ (in the assumption $\sma\Om$). Before the 
construction starts, we prove the principal combinatorial fact.

\ble
\label{comb}
Let\/ $n\in\om$ and\/ $X_u$ be nonempty\/ $\od$ for each\/ 
$u\in 2^n.$ Assume that an\/ $\od$ set\/ $\J{uv}\sq \cD^2$ is 
given for every neighbouring pair of\/ $u,\,v\in 2^n$ so that 
$X_u \J{uv} X_v$.\its
\ben
\def\theenumi{{\rmt\arabic{enumi}.}}
\def\labelenumi{\theenumi}
\item If\/ $u_0\in 2^n$ and\/ $X'\sq X_{u_0}$ is\/ $\od$ and 
nonempty then there exists a system of\/ $\od$ nonempty sets\/ 
$Y_u\sq X_u\;\;(u\in 2^n)$ such that\/ $Y_u \J{uv} Y_v$ holds for 
all neighbouring pairs\/ $u,\,v,$ and in addition $Y_{u_0}=X'$.\its

\item
Suppose that\/ $u_0,\,v_0\in 2^n$ is a neighbouring pair and 
nonempty\/ $\hspace{-1pt}\od\hspace{-1pt}$ 
sets\/ ${X'\sq X_{u_0}}$ and $X''\sq X_{v_0}$ satisfy\/ 
${X' \J{u_0v_0} X''}.$ Then there exists a system of\/ $\od$ 
nonempty sets\/ ${Y_u\sq X_u}$ $(u\in 2^n)$ such that\/ 
${Y_u \J{uv} Y_v}$ holds for all neighbouring pairs\/ $u,v,$ and in 
addition\/ $Y_{u_0}=X',\,\;Y_{v_0}=X''.$
\een
\ele
\proof Notice that 1 follows from 2. Indeed take arbitrary $v_0$ 
such that either $\ang{u_0,v_0}$ or $\ang{v_0,u_0}$ is neighbouring, 
and put respectively 
${X''=\ans{y\in X_{v_0}: \exists\,x\in X'\;(x \J{u_0v_0} y)}},$ or 
${X''=\ans{y\in X_{v_0}: \exists\,x\in X'\;(y \J{v_0u_0} x)}}$. 

To prove item 2, we use induction on $n.$ 

For $n=1$ --- then $u_0=\ang{0}$ and $v_0=\ang{1}$ --- 
we take $Y_{u_0}=Y'$ and $Y_{v_0}=Y''$.

The step. We prove the lemma for $n+1$ provided it has been proved 
for $n;\,\,n\>1.$ The principal idea is to divide $2^{n+1}$ on two 
copies of $2^n,$ minimally connected by neighbouring pairs, and handle 
them more or less separately using the induction hypothesis. The 
two ``copies'' are $U_0=\ans{s\we 0:s\in 2^n}$ and 
$U_1=\ans{s\we 1:s\in 2^n}$. 

The only neighbouring pair that connects $U_0$ and $U_1$ is the pair 
of $\ovu=0^n\we 0$ and $\ovv=0^n\we 1.$ If in fact $u_0=\ovu$ 
and $v_0=\ovv$ then we apply the induction hypothesis (item~1) 
independently for the families of sets ${\ans{X_u:u\in U_0}}$ and 
${\ans{X_u:u\in U_1}}$ and the given sets ${X'\sq X_{u_0}}$ and 
${X''\sq X_{v_0}.}$ Assembling the results, we get nonempty $\od$ 
sets ${Y_u\sq X_u\,\;(u\in 2^{n+1})}$ such that ${Y_u \J{uv} Y_v}$ 
for all neighbouring pairs\/ $u,\,v,$ perhaps with the exception of the 
pair of $u=u_0=\ovu$ and $v=v_0=\ovv,$ and in addition 
$Y_{u_0}=X'$ and $Y_{v_0}=X''.$ 
Thus finally $Y_\ovu \J{\ovu\ovv}Y_\ovv$ by the choice of $X'$ and 
$Y'$. 

It remains to consider the case when both $u_0$ and $v_0$ belong 
to one and the same domain, say to $U_0.$ Then we first apply the 
induction hypothesis (item 2) to the family 
${\ans{X_u:u\in U_0}}$ and the sets ${X'\sq X_{u_0}}$ and 
${X''\sq X_{v_0}.}$ This results in a system of nonempty $\od$ 
sets ${Y_u\sq X_u\;\,(u\in U_0);}$ in particular we get an $\od$ 
nonempty set $Y_\ovu\sq X_\ovu.$ It remains to put 
${Y_\ovv=\ans{y\in X_\ovv:\exists\,x\in Y_\ovu\,
(x\J{\ovu\ovv}y)},}$ so that ${Y_\ovu\J{\ovu\ovv}Y_\ovv,}$ 
and to apply the induction hypothesis (item 1) to the family 
${\ans{X_u:u\in U_1}}$ and the set $Y_\ovv\sq X_\ovv$.\qed

\subsubsection{The construction}

We put $X_\La=X_0.$ 

Now assume that the sets $X_s\,\;(s\in 2^n)$ 
and relations $\J{st}$ for all neighbouring pairs of $s,\,t\in 2^{\<n}$ 
have been defined, and expand the construction at level $n+1.$ 

We first put $A_{s\we i}=X_s$ for all $s\in 2^n$ and 
$i\in\ans{0,1}.$ We also define $\I{uv}=\J{st}$ for any neighbouring 
pair of $u=s\we i,\,\,v=t\we i$ in $2^{n+1}$ other than the pair 
$\ovu=0^n\we 0,\,\,\ovv=0^n\we 1.$ For the latter one (notice 
that $A_{\ovu}=A_{\ovv}=X_{0^n}$) we put $\I{\ovu\ovv}=\oE,$ 
so that $A_u\I{uv} A_v$ holds for all neighbouring pairs of 
$u,\,v\in 2^{n+1}$ including the pair $\ang{\ovu,\ovv}$.

The sets $A_u$ and $\I{uv}$ will be reduced in several steps to 
meet requirements \ref{a}, \ref{b}, \ref{d} and \ref{g1}, \ref{g2}, 
\ref{g3} of Subsection~\ref{embed}.\vom

{\em Part 1}. After $2^{n+1}$ steps of the procedure of 
Lemma~\ref{comb} (item 1) we 
obtain a system of nonempty $\od$ sets 
$B_u\sq A_u\;\,(u\in 2^{n+1})$ such that still $B_u\I{uv} B_v$ 
for all neighbouring pairs $u,\,v$ in $2^{n+1},$ but $B_u\in \cX_{n+1}$ 
for all $u.$ Thus \ref{g1} is fixed.\vom

{\em Part 2}. To fix \ref{g2}, consider an arbitrary pair of 
$u_0=s_0\we 0,$ $v_0=t_0\we 1,$ where $s_0,\,t_0\in 2^n.$ By 
Remark~\ref{newrem} and density of the set $\cP_{n+1}$ there 
exist nonempty $\od$ sets $B'\sq B_{u_0}$ and $B''\sq B_{v_0}$ 
such that ${P=(B'\ti B'')\cap\oE\in \cP_{n+1}}$ and 
$\pri P=B',$ $\prii P=B'',$ so in particular ${B'\oE B''}.$ Now we 
apply Lemma~\ref{comb} (item 1) separately for the two systems 
of sets, 
${\ans{B_{s\we 0}:s\in 2^n}}$ and ${\ans{B_{t\we 1}:t\in 2^n}}$ 
(compare with the proof of Lemma~\ref{comb}~!), and the sets 
$B'\sq B_{s_0\we 0},$ $B''\sq B_{t_0\we 1}$ respectively. 
This results in a system of nonempty 
$\od$ sets ${B'_u\sq B_u}$ ${(u\in 2^{n+1})}$ satisfying 
${B'_{u_0}=B'}$ and ${B'_{v_0}=B'',}$ so that we have 
${(B'_{u_0}\ti B'_{v_0})\cap\oE\in \cP_{n+1},}$ and still 
$B'_u\I{uv} B'_v$ for all neighbouring pairs $u,\,v\in 2^{n+1},$ 
perhaps with the exception of the pair of 
$\ovu=0^n\we 0,\,\,\ovv=0^n\we 1,$ which is the only one that 
connects the two domains. To handle this exceptional pair, 
note that ${B'_{\ovu} \oE B'_{u_0}}$ and ${B'_{\ovv} \oE B'_{v_0}}$ 
(Remark~\ref{newrem} is applied to each of the two domains), 
so that ${B'_\ovu\oE B'_\ovv}$ since ${B'\oE B''}.$ 
We observe that $\I{\ovu\ovv}$ is so far equal to $\oE$. 

After $2^{n+1}$ steps (the number of pairs $u_0,\,v_0$ to be 
considered here) we get a system of nonempty $\od$ sets 
$C_u\sq B_u\;\,(u\in 2^{n+1})$ such that  
$(C_u\ti C_v)\cap\oE\in \cP_{n+1}$ whenever $u(n)\not=v(n),$ 
and still $C_u\I{uv} C_v$ for all neighbouring pairs 
$u,\,v\in 2^{n+1}.$ Thus \ref{g2} is fixed.\vom

{\em Part 3}. We fix \ref{d} for the exceptional neighbouring 
pair of ${\ovu=0^n\we 0},$ ${\ovv=0^n\we 1}.$ Since $\E$ is 
\dd{\top^2}dense in $\oE,$ and ${C_\ovu\oE C_\ovv,}$ the set
${\R{}=(C_\ovu \ti C_\ovv)\cap\E}$ is nonempty. Then some nonempty 
$\od$ set $\I{}\sq \R{}$ is \dd\top separable by Lemma~\ref{dizl}. 
Consider the $\od$ sets $C'=\pri \I{}\,\,(\sq C_\ovu)$ and 
$C''=\prii \I{}\,\,(\sq C_\ovv);$ obviously $C'\I{} C'',$ so that
$C'\I{\ovu\ovv} C''.$ (We recall that at the moment 
$\I{\ovu\ovv}=\oE.$) Using Lemma~\ref{comb} (item 2) again, we 
obtain a system of nonempty $\od$ sets 
$Y_u\sq C_u\;\,(u\in 2^{n+1})$ such that still $Y_u\I{uv} Y_v$ for 
all neighbouring pairs $u,\,v$ in $2^{n+1},$ and $Y_\ovu=C',$  
$Y_\ovv=C''.$ We re--define $\I{\ovu\ovv}$ by 
$\I{\ovu\ovv}=\I{},$ but this keeps 
$Y_\ovu\I{\ovu\ovv} Y_\ovv$.\vom

{\em Part 4}. We fix \ref{g3}. Consider a neighbouring pair $u_0,\,v_0$ 
in $2^{n+1}.$ Then $u_0=0^k\we 0\we r,$ $v_0=0^k\we 1\we r$ for 
some ${k\<n}$ and ${r\in 2^{n-k}}.$ We observe that the temporary 
relation ${\Ip{}=\I{u_0v_0}\cap(Y_{u_0}\ti Y_{v_0})}$ is a nonempty 
(because ${Y_{u_0}\I{u_0v_0} Y_{v_0}}$) $\od$ subset of 
$\J k=\J{0^k\we 0\,,\,0^k\we 1}$ by the construction. Let 
$\I{}\sq \Ip{}$ be a nonempty $\od$ set in $\cX_{n+1}(\J k).$ Now 
put $Y'=\pri \I{}$ and $Y''=\prii\I{}$ (then ${Y'\I{} Y''}$ 
and ${Y'\I{u_0v_0}Y''}$) and 
run Lemma~\ref{comb} (item 2) for the system of sets 
$Y_u\;\,(u\in 2^{n+1})$ and the sets ${Y'\sq Y_{u_0}},$ 
${Y''\sq Y_{v_0}}$. After this define the ``new'' 
$\I{u_0v_0}$ by $\I{u_0v_0}=\I{}$.  

Do this consequtively for all neighbouring pairs; the finally obtained 
sets -- let us denote them by $X_u\,\;(u\in 2^{n+1})$ -- are as 
required. The final relations $\J{uv}\;\,(u,\,v\in 2^{n+1})$ can 
be obtained as the restrictions of the relations 
$\I{uv}$ to $X_u\ti X_v$.\vom

This ends the construction. \vspace{3mm}

This also ends the proof of theorems \ref{mt} and 
\ref{main}.\qed\qed

\end{document}